\documentclass[12pt]{amsart}       
\usepackage{txfonts}
\usepackage{algorithm}
\usepackage{algorithmic}
\usepackage{amssymb}
\usepackage{eucal}
\usepackage{graphicx}
\usepackage{amsmath}
\usepackage{amscd}
\usepackage[all]{xy}           
\usepackage{tikz}
\usepackage{tikz-qtree}
\usepackage{amsfonts,latexsym}
\usepackage{xspace}
\usepackage{epsfig}
\usepackage{float}
\usepackage{color}
\usepackage{fancybox}
\usepackage{colordvi}
\usepackage{multicol}
\usepackage{colordvi}
\ifpdf
  \usepackage[colorlinks,final,backref=page,hyperindex]{hyperref}
\else
  \usepackage[colorlinks,final,backref=page,hyperindex,hypertex]{hyperref}
\fi
\usepackage[active]{srcltx} 
\usepackage{mathrsfs} 

\newtheorem{theorem}{Theorem}[section]
\newtheorem{lemma}[theorem]{Lemma}

\newtheorem{corollary}[theorem]{Corollary}
\newtheorem{prop-def}{Proposition-Definition}[section]
\newtheorem{coro-def}{Corollary-Definition}[section]
\newtheorem{corollary-def}{Corollary-Definition}[section]

\newtheorem{conjecture}[theorem]{Conjecture}
\newtheorem{problem}{Problem}
\theoremstyle{definition}
\newtheorem{defn}[theorem]{Definition}
\newtheorem{remark}[theorem]{Remark}

\newtheorem{tempex}[theorem]{Example}
\newtheorem{tempexs}[theorem]{Examples}
\newenvironment{exam}{\begin{tempex}\rm}{\end{tempex}}

\newcommand{\nc}{\newcommand}
\nc{\tred}[1]{\textcolor{red}{#1}}
\nc{\tblue}[1]{\textcolor{blue}{#1}}
\nc{\tgreen}[1]{\textcolor{green}{#1}}
\nc{\tpurple}[1]{\textcolor{purple}{#1}}
\nc{\btred}[1]{\textcolor{red}{\bf #1}}
\nc{\btblue}[1]{\textcolor{blue}{\bf #1}}
\nc{\btgreen}[1]{\textcolor{green}{\bf #1}}
\nc{\btpurple}[1]{\textcolor{purple}{\bf #1}}
\nc{\NN}{{\mathbb N}}
\nc{\ncsha}{{\mbox{\cyr X}^{\mathrm NC}}} \nc{\ncshao}{{\mbox{\cyr
X}^{\mathrm NC}_0}}

\renewcommand{\frak}{\mathfrak}

\newcommand{\efootnote}[1]{}

\renewcommand{\textbf}[1]{}

\newcommand{\delete}[1]{}

\nc{\mlabel}[1]{\label{#1}}  
\nc{\mcite}[1]{\cite{#1}}  
\nc{\mref}[1]{\ref{#1}}  
\nc{\meqref}[1]{\eqref{#1}}  
\nc{\mbibitem}[1]{\bibitem{#1}} 

\delete{
\nc{\mlabel}[1]{\label{#1}  
{\hfill \hspace{1cm}{\small\tt{{\ }\hfill(#1)}}}}
\nc{\mcite}[1]{\cite{#1}{\small{\tt{{\ }(#1)}}}}  
\nc{\mref}[1]{\ref{#1}{{\tt{{\ }(#1)}}}}  
\nc{\meqref}[1]{\eqref{#1}{{\tt{{\ }(#1)}}}}  
\nc{\mbibitem}[1]{\bibitem[\bf #1]{#1}} 
}

\nc\ns{ns\xspace}

\nc{\opo}{\mathscr{O}}
\nc{\opp}{\mathscr{P}}
\nc{\opf}{\mathscr{F}}
\nc{\optopd}{\mathscr{O}}

\nc{\name}[1]{{\bf #1}}

\nc{\tforall}{\quad \text{for all }}

\nc{\mo}{\calm}
\nc{\barot}{\overline{\otimes}}
\nc{\dm}{\diamond_\mo}
\nc{\oot}{\overline{\ot}}
\nc{\opa}{\ast} \nc{\opb}{\odot} \nc{\op}{\bullet} \nc{\pa}{\frakL}
\nc{\arr}{\rightarrow} \nc{\lu}[1]{(#1)} \nc{\mult}{\mrm{mult}}
\nc{\diff}{\mathfrak{Diff}}
\nc{\opc}{\sharp}\nc{\opd}{\natural}
\nc{\ope}{\circ}
\nc{\dpt}{\mathrm{d}}
\nc{\GS}{Gr\"obner-Shirshov\xspace}
\nc{\gsb}{Gr\"obner-Shirshov basis\xspace}
\nc{\gsbs}{Gr\"{o}bner-Shirshov bases\xspace}
\nc{\diam}{alternating\xspace}
\nc{\Diam}{Alternating\xspace}
\nc{\cdiam}{canonical alternating\xspace}
\nc{\Cdiam}{Canonical alternating\xspace}
\nc{\AW}{\mathcal{A}}
\nc{\mrbo}{modified RBO\xspace }
\nc{\ari}{\mathrm{ar}}
\nc{\lef}{\mathrm{lef}}
\nc{\Sh}{\mathrm{ST}}
\nc{\Cr}{\mathrm{Cr}}
\nc{\st}{{Schr\"oder tree}\xspace}
\nc{\sts}{{Schr\"oder trees}\xspace}
\nc{\vertset}{\Omega} 
\nc{\pb}{{\mathrm{pb}}}
\nc{\Lf}{{\mathrm{Lf}}}
\nc{\lft}{{left tree}\xspace}
\nc{\lfts}{{left trees}\xspace}
\nc{\fat}{{fundamental averaging tree}\xspace}
\nc{\fats}{{fundamental averaging trees}\xspace}
\nc{\avt}{\mathrm{Avt}}
\nc{\rass}{{\mathit{RAss}}}
\nc{\aass}{{\mathit{AAss}}}

\nc{\vin}{{\mathrm Vin}}    
\nc{\lin}{{\mathrm Lin}}    
\nc{\inv}{\mathrm{I}n}
\nc{\gensp}{V} 
\nc{\genbas}{\mathcal{V}} 
\nc{\bvp}{V_P}     
\nc{\gop}{{\,\omega\,}}     

\nc{\bin}[2]{ (_{\stackrel{\scs{#1}}{\scs{#2}}})}  
\nc{\binc}[2]{ \left (\!\! \begin{array}{c} \scs{#1}\\
    \scs{#2} \end{array}\!\! \right )}  
\nc{\bincc}[2]{  \left ( {\scs{#1} \genfrac
    \vspace{-1cm}\scs{#2}} \right )}  
\nc{\bs}{\bar{S}} \nc{\cosum}{\sqsubset} \nc{\la}{\longrightarrow}
\nc{\rar}{\rightarrow} \nc{\dar}{\downarrow} \nc{\dprod}{**}
\nc{\dap}[1]{\downarrow \rlap{$\scriptstyle{#1}$}}
\nc{\md}{\mathrm{dth}} \nc{\uap}[1]{\uparrow
\rlap{$\scriptstyle{#1}$}} \nc{\defeq}{\stackrel{\rm def}{=}}
\nc{\disp}[1]{\displaystyle{#1}} \nc{\dotcup}{\
\displaystyle{\bigcup^\bullet}\ } \nc{\gzeta}{\bar{\zeta}}
\nc{\hcm}{\ \hat{,}\ } \nc{\hts}{\hat{\otimes}}
\nc{\free}[1]{\bar{#1}}
\nc{\uni}[1]{\tilde{#1}} \nc{\hcirc}{\hat{\circ}} \nc{\lleft}{[}
\nc{\lright}{]} \nc{\lc}{\lfloor} \nc{\rc}{\rfloor}
\nc{\curlyl}{\left \{ \begin{array}{c} {} \\ {} \end{array}
    \right .  \!\!\!\!\!\!\!}
\nc{\curlyr}{ \!\!\!\!\!\!\!
    \left . \begin{array}{c} {} \\ {} \end{array}
    \right \} }
\nc{\longmid}{\left | \begin{array}{c} {} \\ {} \end{array}
    \right . \!\!\!\!\!\!\!}
\nc{\onetree}{\bullet} \nc{\ora}[1]{\stackrel{#1}{\rar}}
\nc{\ola}[1]{\stackrel{#1}{\la}}
\nc{\ot}{\otimes} \nc{\mot}{{{\boxtimes\,}}}
\nc{\otm}{\overline{\boxtimes}} \nc{\sprod}{\bullet}
\nc{\scs}[1]{\scriptstyle{#1}} \nc{\mrm}[1]{{\rm #1}}
\nc{\margin}[1]{\marginpar{\rm #1}}   
\nc{\dirlim}{\displaystyle{\lim_{\longrightarrow}}\,}
\nc{\invlim}{\displaystyle{\lim_{\longleftarrow}}\,}
\nc{\mvp}{\vspace{0.3cm}} \nc{\tk}{^{(k)}} \nc{\tp}{^\prime}
\nc{\ttp}{^{\prime\prime}} \nc{\svp}{\vspace{2cm}}
\nc{\vp}{\vspace{8cm}} \nc{\proofbegin}{\noindent{\bf Proof: }}
\nc{\proofend}{$\blacksquare$ \vspace{0.3cm}}
\nc{\modg}[1]{\!<\!\!{#1}\!\!>}
\nc{\intg}[1]{F_C(#1)} \nc{\lmodg}{\!
<\!\!} \nc{\rmodg}{\!\!>\!}
\nc{\cpi}{\widehat{\Pi}}
\nc{\sha}{{\mbox{\cyr X}}}  
\nc{\ssha}{{\mbox{\cyrs X}}} 
\nc{\shpr}{\diamond}    
\nc{\shp}{\ast} \nc{\shplus}{\shpr^+}
\nc{\shprc}{\shpr_c}    
\nc{\msh}{\ast} \nc{\zprod}{m_0} \nc{\oprod}{m_1}
\nc{\vep}{\varepsilon} \nc{\labs}{\mid\!} \nc{\rabs}{\!\mid}
\nc{\sqmon}[1]{\langle #1\rangle}

\nc{\mmbox}[1]{\mbox{\ #1\ }} \nc{\dep}{\mrm{dep}} \nc{\fp}{\mrm{FP}}
\nc{\rchar}{\mrm{char}} \nc{\End}{\mrm{End}} \nc{\Fil}{\mrm{Fil}}
\nc{\Mor}{Mor\xspace} \nc{\gmzvs}{gMZV\xspace}
\nc{\gmzv}{gMZV\xspace} \nc{\mzv}{MZV\xspace}
\nc{\mzvs}{MZVs\xspace} \nc{\Hom}{\mrm{Hom}} \nc{\id}{\mrm{id}}\nc\Idl{{\Id_{\rm Lie}}}
\nc{\im}{\mrm{im}} \nc{\incl}{\mrm{incl}} \nc{\map}{\mrm{Map}}
\nc{\mchar}{\rm char} \nc{\nz}{\rm NZ} \nc{\supp}{\rm Supp}

\nc{\Alg}{\mathbf{Alg}} \nc{\Bax}{\mathbf{Bax}} \nc{\bff}{\mathbf f}
\nc{\bfk}{{\bf k}} \nc{\bfone}{{\bf 1}} \nc{\bfx}{\mathbf x}
\nc{\bfy}{\mathbf y}
\nc{\base}[1]{\bfone^{\otimes ({#1}+1)}} 
\nc{\Cat}{\mathbf{Cat}}

\nc{\detail}{\marginpar{\bf More detail}
    \noindent{\bf Need more detail!}
    \svp}
\nc{\Int}{\mathbf{Int}} \nc{\Mon}{\mathbf{Mon}}
\nc{\rbtm}{{shuffle }} \nc{\rbto}{{Rota-Baxter }}
\nc{\remarks}{\noindent{\bf Remarks: }} \nc{\Rings}{\mathbf{Rings}}
\nc{\Sets}{\mathbf{Sets}} \nc{\wtot}{\widetilde{\odot}}
\nc{\wast}{\widetilde{\ast}} \nc{\bodot}{\bar{\odot}}
\nc{\bast}{\bar{\ast}} \nc{\hodot}[1]{\odot^{#1}}
\nc{\hast}[1]{\ast^{#1}} \nc{\mal}{\mathcal{O}}
\nc{\tet}{\tilde{\ast}} \nc{\teot}{\tilde{\odot}}
\nc{\oex}{\overline{x}} \nc{\oey}{\overline{y}}
\nc{\oez}{\overline{z}} \nc{\oef}{\overline{f}}
\nc{\oea}{\overline{a}} \nc{\oeb}{\overline{b}}
\nc{\weast}[1]{\widetilde{\ast}^{#1}}
\nc{\weodot}[1]{\widetilde{\odot}^{#1}} \nc{\hstar}[1]{\star^{#1}}
\nc{\lae}{\langle} \nc{\rae}{\rangle}
\nc{\lf}{\lfloor}
\nc{\rf}{\rfloor}

\nc{\QQ}{{\mathbb Q}}
\nc{\RR}{{\mathbb R}} \nc{\ZZ}{{\mathbb Z}}

\nc{\cala}{{\mathcal A}} \nc{\calb}{{\mathcal B}}
\nc{\calc}{{\mathcal C}}
\nc{\cald}{{\mathcal D}} \nc{\cale}{{\mathcal E}}
\nc{\calf}{{\mathcal F}} \nc{\calg}{{\mathcal G}}
\nc{\calh}{{\mathcal H}} \nc{\cali}{{\mathcal I}}
\nc{\call}{{\mathcal L}} \nc{\calm}{{\mathcal M}}
\nc{\caln}{{\mathcal N}} \nc{\calo}{{\mathcal O}}
\nc{\calp}{{\mathcal P}} \nc{\calr}{{\mathcal R}}
\nc{\cals}{{\mathcal S}} \nc{\calt}{{\mathcal T}}
\nc{\calu}{{\mathcal U}} \nc{\calw}{{\mathcal W}} \nc{\calk}{{\mathcal K}}
\nc{\calx}{{\mathcal X}} \nc{\CA}{\mathcal{A}}

\nc{\fraka}{{\mathfrak a}} \nc{\frakA}{{\mathfrak A}}
\nc{\frakb}{{\mathfrak b}} \nc{\frakB}{{\mathfrak B}}
\nc{\frakc}{{\mathfrak c}}
\nc{\frakD}{{\mathfrak D}} \nc{\frakF}{\mathfrak{F}}
\nc{\frakf}{{\mathfrak f}} \nc{\frakg}{{\mathfrak g}}
\nc{\frakH}{{\mathfrak H}} \nc{\frakI}{{\mathfrak I}}
\nc{\frakL}{{\mathfrak L}}
\nc{\frakM}{{\mathfrak M}} \nc{\bfrakM}{\overline{\frakM}}
\nc{\frakm}{{\mathfrak m}} \nc{\frakP}{{\mathfrak P}}
\nc{\frakN}{{\mathfrak N}} \nc{\frakp}{{\mathfrak p}}
\nc{\frakS}{{\mathfrak S}} \nc{\frakT}{\mathfrak{T}}

\nc{\BS}{\mathbb{S
}}

\font\cyr=wncyr10 \font\cyrs=wncyr7
\nc{\li}[1]{\textcolor{red}{#1}}
\nc{\lir}[1]{\textcolor{red}{Li: #1}}
\nc{\xing}[1]{\textcolor{blue}{Xing: #1}}
\nc{\hu}[1]{\textcolor{purple}{Huhu: #1}}
\nc{\Hu}[1]{\textcolor{purple}{#1}}
\nc{\UN}{U_{N}}
\nc{\FN}{F_{\mathrm M}}
\nc{\altx}{\Lambda}
\nc{\spr}{\cdot}
\nc{\rts}{\stackrel{\rightarrow}{\shpr}}
\nc{\ox}{\overline{\frak x}}
\nc{\oX}{\overline{X}}

\nc{\tree}{\mathbb{T}} \nc{\trees}{\mathbb{T}^\star} \nc{\treehh}{\bfktree^\star}
\nc{\treess}{\mathbb{T}^{\star_1, \star_2}}\nc{\etree}{\mathbb{1}}\nc{\bfktree}{\mathscr{T}}
\nc{\lbar}[1]{\overline{#1}} \nc{\sub}[1]{[#1]}\nc{\suba}[1]{|_{#1}} \nc{\lea}{{\rm L}}\nc{\fg}{{\rm fg}}\nc{\wt}{{\rm wt}} \nc{\degr}{{\rm deg}} \nc{\re}[1]{R(#1)}
\nc{\oid}{\mathrm{OId}} \nc{\irr}{{\rm Irr}}\nc{\irrl}{{\rm Irr_{Lie}}}   \nc{\pis}{\Pi_S} \nc{\dps}{\dotplus}
\nc{\astarrow}{\overset{\raisebox{-2pt}{{\scriptsize $\ast$}}}{\rightarrow}}
\nc{\tvarrow}[3]{#1 \overset{(t,v)}{\longrightarrow}_{#3} #2}\nc{\gs}{Gr\"{o}bner-Shirshov\xspace}
\nc{\rcp}{Rota's Program on Algebraic Operators\xspace}
\nc{\rcpo}{Rota's Program on Algebraic Operators for operads\xspace}
\nc\olie{operated Lie algebra\xspace}
\nc\olies{operated Lie algebras\xspace}
\nc\Olies{Operated Lie algebras\xspace}
\nc{\inp}{\mathrm{In}}\nc{\kirr}{\bfk\irr(S)} \nc{\rbp}{\mathscr{RBP}}
\nc\dlex{<_{dlex}} \nc\plex{<_{{\rm plex}}} \nc\pelex{\leq_{{\rm plex}}} \nc\cond{{j'<j \, \text{or}\atop  i'<i,\,j'=j}}

\nc{\compcs}{compatible compositions }
\nc{\compc}{compatible composition }
\nc\bws[1]{{\lfloor#1\rfloor}}\nc\opmx{\bfk\mathfrak{M}(X)}\nc\opmxm{\mathfrak{M}(X)}\nc\opm[1]{\mathfrak{M}(#1)}
\nc\opmz{\bfk\mathfrak{M}(Z)}\nc\opmzm{\mathfrak{M}(Z)}
\nc\Id{\rm Id}\nc\sopm[1]{\mathfrak{S}^\star(#1)}
\nc\blw[1]{\lfloor#1\rfloor}\nc\plie[1]{\mathfrak{S}(#1)}\nc\plien[1]{\mathcal{N}(#1)}\nc\nas[1]{{#1}^\ast}
\nc\ordc{>_{{\rm Dl}}} \nc\ordqc{\geq_{{\rm Dl}}}
\nc\ord{>_{{\rm DL }}}\nc\ordq{\geq _{\rm DL}}
\nc\ordd{>_{{\rm dt}}}\nc\ordqd{\geq_{{\rm dt}}}
\nc\ordb{>_{{\rm Dl}}}\nc\ordqb{\geq_{{\rm Dl}}}
\nc\alsw[1]{{\rm ALSW}(#1)} \nc\nlsw[1]{{\rm NLSW}(#1)}\nc\alsbw[1]{{\rm ALSBW}_{\ordq}(#1)} \nc\nlsbw[1]{{\rm NLSBW}_{\ordq}(#1)}
\nc\alsbwo[2]{{\rm ALSBW}_{#2}(#1)} \nc\nlsbwo[2]{{\rm NLSBW}_{#2}(#1)}
\nc\oplie{{\rm OLie}(X)}\nc\sopma[1]{\mathfrak{M}^\star(#1)}\nc\Pia[1]{\Pi_{#1}^{\rm ass}}\nc\Pil[1]{\Pi_{#1}^{\rm Lie}}
\nc\opliez{{\rm OLie}(Z)}\nc\opliex{{\rm OLie}(\{x,y\})}
\nc\nlsbwd[1]{\nlsw{\Delta{(#1)}}}\nc\der[2]{{#1}^{(#2)}}

\begin{document}
\title[Rota's program on algebraic operators]{Rota's program on algebraic operators, rewriting systems and Gr\"obner-Shirshov bases}

\author{Xing Gao}
\address{School of Mathematics and Statistics,
Key Laboratory of Applied Mathematics and Complex Systems,
Lanzhou University, Lanzhou, Gansu 730000, P. R. China}
\email{gaoxing@lzu.edu.cn}

\author{Li Guo}\thanks{Corresponding author}
\address{
	Department of Mathematics and Computer Science,
	Rutgers University,
	Newark, NJ 07102, USA}
\email{liguo@rutgers.edu}

\author{Huhu Zhang}
\address{School of Mathematics and Statistics,
	Lanzhou University, Lanzhou, Gansu 730000, P. R. China}
\email{zhanghh17@lzu.edu.cn}

\date{\today}
\begin{abstract}
Many years ago, Rota proposed a program on determining algebraic identities that can be satisfied by linear operators. After an extended period of dormant, progress on this program picked up speed in recent years, thanks to perspectives from operated algebras and Gr\"obner-Shirshov bases. These advances were achieved in a series of papers from special cases to more general situations. These perspectives also indicate that Rota's insight can be manifested very broadly, for other algebraic structures such as Lie algebras, and further in the context of operads. This paper gives a survey on the motivation, early developments and recent advances on Rota's program, for linear operators on associative algebras and Lie algebras. Emphasis will be given to the applications of rewriting systems and Gr\"obner-Shirshov bases. Problems, old and new, are proposed throughout the paper to prompt further developments on \rcp.
\end{abstract}

\subjclass[2010]{
16W99, 
16S15,	
16Z10,	
05C05,   
17B40,	
18M60,	
68Q42	
}

\keywords{\rcp; rewriting systems, Gr\"obner-Shirshov basis; operated associative algebras; operated Lie algebras; operads; differential type operators; Rota-Baxter type operators}
\maketitle

\tableofcontents

\setcounter{section}{0}

\section{Introduction}
This paper provides an introduction, a survey and some questions on Rota's Program on Algebraic Operators.

\subsection{Roles played by linear operators}
Throughout the development of mathematics, linear operators satisfying operator identities have played a pivot role. Such operators can be defined for various algebraic structures. For this Introduction, we will restrict ourselves to the associative context for the sake of historical perspective and simplicity. See Section~\mref{sec:lie} for extending Rota's program to the Lie algebraic context.

We list some of the operators below as motivation for Rota's program. 
\begin{enumerate}
	\item
From the study of algebra and Galois theory arose the homomorphisms of algebras and groups, satisfying the operator identity
\begin{equation*}
	P(xy)=P(x)P(y).
\end{equation*}
\item In the study of analysis, there are the differential operator (derivation), satisfying the Leibniz rule
\begin{equation*}
	d(xy)=d(x)y+xd(y)
\mlabel{eq:leib}
\end{equation*}
and
\item the integral operator, satisfying the integration-by-parts formula
\begin{equation*}
	\int_a^x f'(t)g(t)dt=f(x)g(x)-f(a)g(a) - \int_a^x f(t)g'(t)dt.
 \mlabel{eq:ibp}
\end{equation*}
\item In his study of probability, G. Baxter~\mcite{Ba} introduced the Rota-Baxter operator, satisfying the operator identity
\begin{equation*}
	P(x)P(y)=P(xP(y))+P(P(x)y)+\lambda P(xy).
 \mlabel{eq:rbo}
\end{equation*}
\item
The Reynolds operator has its origin in fluid mechanics~\mcite{Re} and is defined by the operator identity
\begin{equation*}
	R(x)R(y)=R(xR(y))+R(R(x)y)-R(R(x)R(y)).
\mlabel{eq:rey}
\end{equation*}
\item
The notion of averaging operator was explicitly defined by Kamp\'{e}de F\'{e}riet~\mcite{Kf} satisfying
\begin{equation*}
P(x)P(y) = P(xP(y)) = P(P(x)y).
\mlabel{eq:ave}
\end{equation*}
It was already implicitly studied by O. Reynolds~\mcite{Re} in 1895 in turbulence theory under the disguise of a Reynolds operator,
since an idempotent Reynolds operator is an averaging operator.
\item
Studies in mathematics, aerodynamics and signal processing gave rise to the Hilbert transformation (for $\lambda=1$)~\mcite{Co,Tri}  (later called modified Rota-Baxter operator) satisfying
\begin{equation*}
	P(x)P(y)=P(xP(y))+P(P(x)y)+ \lambda xy,
\end{equation*}
which is also called a convolution theorem in~\mcite{Tri}.
\end{enumerate}
The modified Rota-Baxter operator is related to a Rota-Baxter operator by a linear transformation.
Its recent name modified Rota-Baxter operator was adapted from the modified Yang-Baxter equation in Lie algebra~\mcite{STS}, stemmed from integrable systems.

\subsection{\rcp}
All these operators were known to Gian-Carlo Rota, thanks to his broad interests spanning from functional analysis, probability, lattice theory, algebra and foremost combinatorics.
While other researchers had taken operator identities for granted and moved on with their studies of particular interests, Rota looked at such identities from a more transcendental
perspective, probably because of his expertise as a Professor of Philosopher as well as one of Applied Mathematics at MIT. Indeed, according to~\mcite{KY}, a motivation for Rota's change of focus from his earlier research interest in analysis was his recognition of the central role played by the operator identities in analysis and probability.

Rota's program on operator identities was formulated in his 1995 survey paper~\cite{Ro2}:
\begin{quote}
	In a series of papers, I have tried to show that other linear operators satisfying algebraic identities may be of equal importance in studying certain algebraic phenomena, and I have posed the problem of finding all possible algebraic identities that can be satisfied by a linear operator on an algebra.
\end{quote}

Thus the Rota's program aims at not only a summary of known algebraic operator identities, but all the potential algebraic operator identities that might arise in mathematical research.

To fix the terminology, we call a linear operator \name{algebraic} if the operator satisfies an algebraic operator identity. See Definition~\mref{de:algop} for the precise definition. Then we can call Rota's program stated above the \name{\rcp}.

Since Rota proposed his program, further exciting applications of the linear operators that Rota noted above have been found. Differential algebra, originated from an algebraic study of differential equation by Ritt~\mcite{Ri} in the 1930's,  have been developed by the school of Kolchin and many others into a vast area of research~\mcite{CGKS,Ko,SP}  with applications ranging from logic to arithmetic geometry and mechanical proof of geometric theorems~\mcite{Wu,Wu2}.
After the promotion of Rota~\mcite{Ro2}, Rota-Baxter algebra experienced a remarkable renaucence this century, most notably by its applications to renormalization of quantum field theory and Stotisctic process~\mcite{BHZ,CK}. See~\mcite{Gub} for an introduction.
See~\mcite{ZGG} and the references therein for later progresses on Reynolds algebras. Likewise see~\mcite{PG} for average algebras and see~\mcite{ZGG1,ZGG2} for modified Rota-Baxter operators.

At the same time, new linear operators have merged. They include but not limited to
\begin{enumerate}
	\item
	differential operators with weights~\mcite{GK}, defined by
\begin{equation*}
	d(xy)=d(x)y+xd(y)+\lambda d(x)d(y),
\mlabel{eq:diffw}
\end{equation*}
putting the derivation and difference operator in the same framework;
	\item Nijenhuis operators, satisfying the operator identity
\begin{equation*}
P(x)P(y) = P(xP(y) + P(x)y - P(xy)),
\mlabel{eq:nij}
\end{equation*}
introduced by Carin\~{e}na {\em et al.}~\mcite{CGM} to study quantum bi-Hamiltonian systems and constructed by analogy with Poisson-Nijenhuis geometry, from the relative Rota-Baxter algebras~\mcite{Uc};
\item TD operators, characterized by the operator identity
\begin{equation*}
P(x)P(y) = P(xP(y) + P(x)y - xP(1)y),
\mlabel{eq:td}
\end{equation*}
introduced by P.~Lerous in a combinatorial study~\mcite{Le}.
\end{enumerate}

Such developments further showed that operator identities are ubiquitous in mathematics and motivated us to revisit \rcp.
In this process, we found that the understanding of \rcp depends on two
key points. First, we
need a suitable framework to formulate precisely what is an
``operator identity,'' and second, we need to determine key
properties that characterize the classes of operator identities
that are of interest to other areas of mathematics, such as those
listed above.

For the first point, we noted that a simplified but analogous
framework had already been formulated in the 1960s and subsequently
explored with great success. This is the study of PI-rings and
PI-algebras, whose elements satisfy a set of polynomial identities,
or PIs for short ~\mcite{Pr,Ro,DF}.

With the extra action of linear operators, the suitable structure is the operated algebras, defined to be an (associative) algebra equipped with a linear operator free of other restraints. Then the free object can be realized as the operated polynomial algebras (OPIs)~\mcite{Guop}, consisting of polynomials with brackets, called operated polynomials.

Thus operator identities should be taken to be OPIs in
the operated polynomial algebra. Then the second point mentioned above may be
interpreted as follows: among all OPIs, which ones are particularly
consistent with the associative algebra structure so that they were singled out for study?
The general idea is that the operator identities should be consistent with the defining conditions of the algebraic structure, in our case the law of associativity.
We realize this idea by making use of two related theories: rewriting
systems and Gr\"{o}bner-Shirshov bases.

As testing grounds for the general framework and methods, two classes of operators were investigated. The first class is for the differential type operators~\cite{GSZ} and the second class is for the Rota-Baxter type operators~\cite{ZGGS}. Through the investigation of these two classes, it became evident that \rcp is intimately related to convergent rewriting systems and Gr\"obner-Shirshov bases. This connection was confirmed by the more general linear operators~\cite{GG}, where it was shown that for a given family of operator identities, the property that they give a convergent rewriting system is equivalent to the property that they form a Gr\"obner-Shirshov basis of the corresponding operated ideal.

Up to this point everything is considered for associative algebras, since \rcp is for linear operators defined on associative algebras.
But important linear operators have also been defined for other algebraic structures, in particular Lie algebras. For example, derivations on Lie algebras were fundamental for their complexes and homologies. Further, the operator form of the classical Yang-Baxter equation can be given by Rota-Baxter operators and more generally $\mathcal O$-operators on Lie algebras~\mcite{BGN,Bo,Kup,STS,Uc}. Thus it is natural to ask the same questions for other types of algebras, with associativity replaced by the Jacobi identity in the case of Lie algebras.
Then a natural question is a suitable analog of \rcp for other kinds of algebras. This motivates us to put \rcp in the framework of Lie algebras and the perspective of operads.

\subsection{Layout of the paper}
The purpose of this paper is to give an introduction and survey of the recent progresses on \rcp~\mcite{GG,GSZ,ZGGS}, as well as the earlier work of Freeman~\mcite{Fr} from another point of view. We also propose several problems for further research in this promising direction.

The paper is organized as follows.

In Section~\mref{sec:rcpa}, we review the structure of operated algebra in which to consider the operator identities in \rcp. Thus a generic operator identity is an element in the free operated polynomial algebra. We then introduce the methods of rewriting systems and Gr\"obner-Shirshov bases as tools to identify operator identities to be determined in \rcp.
We then apply the general framework and methods to two families of operator identities. One is the differential type operator identities and one is the Rota-Baxter type operators identities. In each case, we introduce the notations, give the equivalence of the desired operator identities with the convergence of a suitable rewriting system and with the existence of a Gr\"obner-Shirshov basis. We also give a conjectured list of such operator identities and provide a uniform construction of the free objects.
We then consider \rcp for more general operator identities and establish a similar equivalence among operator identities that define convergent rewriting systems and those that possess Gr\"obner-Shirshov bases.

In Section~\mref{sec:lie}, we outline a Lie algebraic theory for \rcp. As it turns out, a largely parallel, though more complex, framework exists.

As an appendix, in Section~\mref{sec:free}, we summarize the work of Freeman~\mcite{Fr} on a classification of some of the classical operator identities by means of the associative products on their graphs.

The operadic approach to \rcp is out of the scope of this paper and will be presented elsewhere.

\smallskip

\noindent
{\bf Notation.}
Throughout this paper, let $\bfk$ be a field, which will be the base field of all algebras, tensor products, as well as linear maps.
By an algebra we mean a unitary associative $\bfk$-algebra.
For a set $X$, we use $\bfk X$ to denote the vector space spanned by $X$.

\section{\rcp for associative algebras}
\mlabel{sec:rcpa}
In this section we summarize the progresses made in recent years in formulating and resolving the problem posed by Rota in his program.

\subsection{Operated algebras}\mlabel{ssec:oalg}
In order to make progress on Rota's Program on Algebraic Operators, we first need to make precise the meaning of the program. In particular, we need a suitable framework to formulate precisely
what is an operator identity. A prototype can be found in an algebraic identity satisfied by an algebra, which is defined to be a polynomial identity in a noncommutative polynomial algebra, as a realization of a free (associative) algebra.
However, an algebraic identity in \rcp involves a linear operator. Thus we take an algebraic identity
satisfied by an operator to be an element in a free object in the category of algebras with an operator, or operated algebras. We start with the basic definitions.

\begin{defn}~\mcite{GG,Guop,Gub,Ku}
\begin{enumerate}
\item An {\bf operated monoid (resp. operated algebra)} is a monoid (resp. algebra) $U$ together with a map (resp. linear map)
$P_U: U\to U$.
\item A morphism from an operated monoid (resp. algebra) $(U,P_U)$ to an operated monoid (resp. algebra) $(V,P_V)$
is a monoid (resp. algebra) homomorphism $f : U\to V $ such that
$f\circ P_U=P_V\circ f.$
\end{enumerate}
\end{defn}

Now we review the construction of free operated algebras in terms of bracketed words.
First for any set $Y$, let  $\bws{Y}$  be the set $\{\bws{y}\, | \,y \in Y \}$ ,
which is just another copy of $Y$ whose elements are denoted by $\bws{y}$ for distinction.

Now let $X$ be a set and $M(X)$ the free monoid on $X$ with identity 1.
The free operated monoid over $X$ can be naturally constructed by the limit of a directed system
$$\Big\{\iota_n: \mathfrak{M}_n\to \mathfrak{M}_{n+1}\Big\}_{n=0}^\infty$$
of free monoids $\mathfrak{M}_n$, where the transition morphisms $\iota_n$ will be natural embeddings. For this purpose,
let $\mathfrak{M}_0=M(X)$ and let
$$\mathfrak{M}_1:=M(X\cup\bws{\mathfrak{M}_0}).$$
Let $\iota_0$ be the natural embedding $\iota_0:\mathfrak{M}_0\hookrightarrow \mathfrak{M}_1$ from the inclusion $X\hookrightarrow X\cup \lc \frakM_0\rc$.  Assuming by induction that, for any given $n\geq0$, we have defined,
for $0 \leq i \leq n + 1$, the free monoids $\mathfrak{M}_i$ with the properties that for $0 \leq i \leq n$, we have $\mathfrak{M}_{i + 1} = M ( X\cup \bws{\mathfrak{M}_i })$ and natural embeddings
$\iota_i : \mathfrak{M}_i \to \mathfrak{M}_{i + 1}.$
Then let
$$\mathfrak{M}_{n+2}:=M(X\cup\bws{\mathfrak{M}_{n+1}}).$$
The identity map on $X$ and the embedding $\iota_n$ together induce an injection
$$\iota_{n+1}:X\cup\bws{\mathfrak{M}_{n}}\hookrightarrow X\cup\bws{\mathfrak{M}_{n+1}},$$
which, by the functoriality of taking free monoids, extends to an embedding (still denoted by $\iota_{n+1}$) of free monoids
$$\iota_{n+1}:\mathfrak{M}_{n+1}=M(X\cup\bws{\mathfrak{M}_{n}})\hookrightarrow M( X\cup\bws{\mathfrak{M}_{n+1}})= \mathfrak{M}_{n+2}.$$
This completes the construction of the directed system.
Finally we define the monoid
$$\mathfrak{M}(X):=\bigcup_{n\geq0}\mathfrak{M}_{n},$$
whose elements are called {\bf bracketed words} or {\bf bracketed monomials on $X$}.

A non-unit element $w$ of $\opmxm$ can be uniquely expressed in the form
$$w = w_1\ldots w_k$$
for some $k$ and some $w_i\in X\cup \bws{\opmxm}$, for $ 1 \leq i \leq k$.
In this case, we call $w_i$ {\bf prime} and $k$ the {\bf breadth} of $w$, denoted by $|w|$. If $w = 1$, we define $|w| = 0$.

\begin{lemma}\mcite{Guop}
Let $i:X\hookrightarrow \mathfrak{M}(X)$ and $j:\mathfrak{M}(X)\hookrightarrow \bfk\mathfrak{M}(X)$ be the natural embeddings. Then
\begin{enumerate}
  \item the triple $(\mathfrak{M}(X), P:= \lc\, \rc, i)$  is the free operated monoid on $X$; and
  \item the triple $(\bfk\mathfrak{M}(X), P, j\circ i)$ is the free operated algebra on $X$, where $P$ is the linear operator induced by $\bws{\,}$. \mlabel{item:freeas}
\end{enumerate}
\mlabel{lem:freeas}
\end{lemma}

Thanks to the above construction of free operated algebras, we now can characterize precisely what an operator identity is in \rcp.

Let $\phi=\phi(x_1,  \ldots ,x_k)\in\bfk \mathfrak{M}(X)$. We call $\phi=0$ (or simply $\phi$) an {\bf operated polynomial identity (OPI)}.
By the universal property of the free operated algebra $\bfk \mathfrak{M}(X)$, for any operated algebra $(R, P)$ and any map
$$\theta : \{x_1,  \ldots ,x_k\}\to R,\quad x_i\mapsto r_i, i=1,\ldots,k,$$
there is a unique morphism
$$\tilde{\theta}:\bfk\mathfrak{M}(X)\to R$$
of operated algebras that extends the map $\theta$.
By a slight abuse of notation, define the {\bf evaluation map}
$$\phi(r_1,  \ldots ,r_k):=\tilde{\theta}(\phi(x_1,  \ldots ,x_k)).$$

\begin{defn}
 Let $\phi\in \opmx$ and $(R,P)$ be an operated algebra. We say that
  $R$ is a {\bf $\phi$-algebra} and $P$ is a {\bf $\phi$-operator}, if
  $$\phi(r_1, \ldots ,r_k)=0,\quad \text{for all } r_1, \ldots ,r_k\in R.$$
  More generally, for a subset $\Phi\subseteq\opmx$, we call $R$ (resp. $P$) a {\bf $\Phi$-algebra} (resp. {\bf $\Phi$-operator}) if $R$ (resp. $P$)
is a $\phi$-algebra (resp. $\phi$-operator) for each $\phi\in\Phi$.
\end{defn}

\begin{defn}
Let $R$ be an algebra. A linear operator $P$ on $R$ is called {\bf algebraic} if there is $0\neq \Phi\subseteq \opmx$ such that $P$ is a $\Phi$-operator.
\mlabel{de:algop}
\end{defn}

See~\mcite{GGL} for a related discussion in the context of integral equations.

\begin{exam}
All the operators in the introduction are algebraic.
\end{exam}

\subsection{Rewriting systems and \gsbs}
Now we review two general theories that will be applied to study \rcp.
\subsubsection{Rewriting systems}
Let us first recall some basic concepts of rewriting systems from~\mcite{BN,ZGGS,GG}.

\begin{defn}
Let $V$ be a vector space with a linear basis $W$.
\begin{enumerate}
\item For $f=\sum\limits_{w\in W}c_w w \in V$ with $c_w\in \bfk$, the {\bf  support} $\supp(f)$ of $f$ is the set $\{w\in W\,|\,c_w\neq 0\}$. By convention, we take $\supp(0) = \emptyset$.
\item Let $f, g\in V$. We use $f \dps g$ to indicate the property that $\supp(f) \cap \supp(g) = \emptyset$. If this is the case, we say $f \dps g$ is a {\bf  direct sum} of $f$ and $g$ and use $f\dps g$ to denote the sum $f+ g$.

\item For $f \in V$ and $w \in \supp(f)$ with the coefficient $\alpha_w$, write $R_w(f) := \alpha_w w -f \in V$. So $f = \alpha_w w \dps (-R_w(f))$.
\end{enumerate}
\mlabel{def:dps}
\end{defn}

\begin{defn}
Let $V$ be a vector space with a linear basis $W$.
\begin{enumerate}
\item  A {\bf  term-rewriting system $\Pi$ on $V$ with respect to $W$} is a binary relation $\Pi \subseteq W \times V$. An element $(t,v)\in \Pi$ is called a (term-) rewriting rule of $\Pi$, denoted by $t\to v$.

\item The term-rewriting system $\Pi$ is called {\bf simple with respect to $W$} if $t \dps v$ for all $t\to v\in \Pi$.

\item If $f = \alpha_t t\dps (-R_t(f))\in V$, using the rewriting rule $t\to v$, we get a new element $g:= c_t v - R_t(f) \in V$, called a {\bf one-step rewriting} of $f$ and denoted $f \to_\Pi g$ or $\tvarrow{f}{g}{\Pi}$.  \label{item:Trule}

\item The reflexive-transitive closure of $\rightarrow_\Pi$ (as a binary relation on $V$) is denoted by $\astarrow_\Pi$ and, if $f \astarrow_\Pi g$, we say {\bf  $f$ rewrites to $g$ with respect to $\Pi$}. \label{item:rtcl}

\item Two elements  $f, g \in V$ are {\bf  joinable} if there exists $h \in V$ such that $f \astarrow_\Pi h$ and $g \astarrow_\Pi h$; we denote this by $f \downarrow_\Pi g$.

\item A {\bf fork}  is a pair of distinct reduction sequences $(f \astarrow_\Pi g_1, f \astarrow_\Pi g_2)$  starting from the same  $f\in V$. The fork is called {\bf joinable} if $g_1 \downarrow_\Pi g_2$.

\item An element $f\in V$ is {\bf a normal form} if no more rules from $\Pi$ can apply.
\end{enumerate}
\label{def:ARSbasics}
\end{defn}

\begin{defn}
A term-rewriting system $\Pi$ on $V$ is called
 \begin{enumerate}
 \item {\bf  terminating} if there is no infinite chain of one-step rewritings
 $$f_0 \rightarrow_\Pi f_1 \rightarrow_\Pi f_2 \rightarrow_\Pi \cdots \quad.$$
\item {\bf  confluent} if every fork  is joinable.
\item {\bf  convergent} if it is both terminating and confluent.
\end{enumerate}
\label{def:ARS}
\end{defn}

\subsubsection{\gsbs}
Next we recall some concepts of \gsbs. See~\cite{BC,GG,GSZ,ZGGS} for further details.

\begin{defn}
Let $Z$ be a set, let $\star$ be a symbol not in $Z$, and let $Z^\star = Z\cup \{\star\}$.
\begin{enumerate}
  \item  By a {\bf $\star$-bracketed word} on $Z$, we mean any word in  $\opm{Z^\star}$ with exactly one occurrence of $\star$. The set of all
$\star$-bracketed words on $Z$ is denoted by $\sopma{Z}$.
  \item  For $q\in \sopma{Z}$ and $u \in\opm{Z}$, we define $q\suba{\star\mapsto u}$ (or $q\suba{u}$) to be the bracketed word on $Z$
obtained by replacing the symbol $\star$ in $q$ by $u$.
  \item  For $q\in \sopma{Z}$ and $s =\sum_i c_i u_i \in\opmz$, where $c_i\in\bfk$ and $u_i\in\opmzm$, we define
$$q\suba{s}:=\sum_i c_i q\suba{u_i}.$$
  \item A bracketed word $u \in\opmzm$ is a {\bf subword} of another bracketed word $ w\in\opmzm$ if
$w = q\suba{u}$ for some $q \in\sopma{Z}$.
\end{enumerate}
\end{defn}

A monomial order is a well order that is compatible with all operations in the algebraic structure. More precisely, we give

\begin{defn}
Let $Z$ be a set. A {\bf monomial order} on $\opmzm$ is a well order $\geq$ on $\opmzm$ such that
\begin{equation*}
\quad  u > v \Rightarrow  q\suba{u} > q\suba{v},\,\text{ for all } u,v\in \opmzm \text{ and } q\in\sopma{Z}.
\end{equation*}
\end{defn}

Given a monomial order $\geq$ and an element $f\in \opmzm$, we let $\lbar{f}$ denote the leading bracketed word (monomial) of $f$.
If the coefficient
of $\lbar{f}$ in $f$ is 1, we call $f$ {\bf monic with respect to the monomial order $\geq$}. A subset $S$ of $\opmz$ is called {\bf monic} if each nonzero element $s$ in $S$ is monic.
We say an element $f\in\opmz$ is in {\bf normal $\phi$-form} if no monomial of $f$ contains any subword of the form $\bar{\phi}$.

We are going to expose the key concept of \gsbs on $\bfk\frakM(Z)$.

\begin{defn}
Let $\geq$ be a monomial order on $\opmzm$ and $f, g \in\opmz$ be monic.
\begin{enumerate}
  \item  If there are $w, u, v\in \opmzm$ such that $w = \lbar{f}u = v\lbar{g}$ with max$\{|\lbar{f}|, |\lbar{g}|\} < |w| <
|\lbar{f}| + |\lbar{g}|$, we call
$$(f,g)^{u,v}_w:= fu - vg$$
the {\bf intersection composition of $f$ and $g$ with respect to $w$}.
\mlabel{item:intcomp}
  \item  If there are $w\in\opmzm$ and $q\in\sopma{Z}$ such that $w = \lbar{f} = q\suba{\lbar{ g}},$ we call
$$(f,g)^q_w := f - q\suba{g}$$
the {\bf including composition of $f$ and $g$ with respect to $w$}.
\mlabel{item:inccomp}
\end{enumerate}
\mlabel{defn:comp}
\end{defn}

\begin{defn}
Let $Z$ be a set and $\geq$ be a monomial order on $\opmzm$. Let $S\subseteq\opmz$ be monic.
\begin{enumerate}
  \item  An element $f\in\opmz$ is called {\bf  trivial modulo $(S, w)$} if
$$f =\sum_ic_i q_i\suba{s_i} \, \text{ with } \lbar{q_i\suba{s_i}} < w, \,\text{ where } c_i\in\bfk, q_i\in\sopma{Z}, s_i\in S.$$
  \item   We call $S$ a {\bf \gsb} in $\opmz$ with respect
to $\geq$ if, for all pairs $f, g \in S$, every intersection composition of the form $(f, g)^{u,v}_w$
is trivial modulo $(S, w)$, and every including composition of the form $(f, g)^q_w$ is trivial modulo $(S, w)$.
\end{enumerate}
\end{defn}

\subsection{Differential type OPIs}
Two particular classes of operators, namely, those that generalize the differential operator or the Rota-Baxter operator were studied in~\mcite{GSZ, ZGGS} respectively.
The first class is the differential type operators.

\begin{defn}\mcite{GSZ}
A {\bf differential type OPI} is
$$\phi(x , y ) := \bws{xy}- N(x , y )\in\bfk\opm{\{x,y\}},$$
where
\begin{enumerate}
  \item $N(x, y )$ is multi-linear in $x$ and $y$ ;
  \mlabel{item:a}
  \item $N(x, y )$ is a normal $\phi$-form, that is, $N(x , y )$  does not contain any subword of the form $\bws{uv}$, for any non-units $u,v\in \bfk\opm{\{x,y\}}$;
  \mlabel{item:b}
  \item For any set $Z$ with $u, v, w \in \opmz\backslash\{1\}$,
$N(uv,w) - N(u,vw)\astarrow_{\Pia{S_\phi}} 0,$
for the term-rewriting system
\begin{equation}
	\Pia{S_\phi}:=\{q\suba{\bws{u  v}}\rightarrow q\suba{N(u  ,v )}\,|\,q\in\sopma{Z},u,v\in\opm{Z} \}.
	\mlabel{eq:difftr0}
\end{equation}
  \mlabel{item:c}
\end{enumerate}
A linear operator satisfying a differential type OPI is called {\bf a differential type operator}.
\mlabel{defn:difftyp}
\end{defn}

\begin{remark}
Condition~\meqref{item:a} is imposed since we are only interested in linear operators. Condition~\meqref{item:b}
is needed to avoid infinite rewriting under $\Pia{S_\phi}$. Condition~\meqref{item:c} is needed so that $\bws{( uv ) w } = \bws{ u ( vw )}$.
\end{remark}

\begin{defn}
The {\bf operator degree} of a monomial in $\opmx$
is the total number that the operator $\bws{\, }$ appears in the monomial. The {\bf operator degree } of a polynomial $f$ in $\opmx$ is the maximum of the
operator degrees of the monomials appearing in $f$.
\end{defn}

The following is a list of differential type operators with operator degrees not exceeding two. See~\mcite{GSZ} for more details.

\begin{theorem} {\em (Classification of differential type operators)}~\mcite{GSZ}
Let $a , b , c , e\in\bfk$. The OPI $\phi(x,y) := \bws{xy} - N(x ,y)\in\bfk\opm{\{x,y\}}$, where $N( x , y )$  is taken from the list below, is of differential type.
\begin{enumerate}
  \item $b ( x \bws{ y }+\bws{ x } y )+ c \bws{ x }\bws{ y }+ exy $ where $b^2 = b + ce$,
  \item  $ce^2 yx + exy + c \bws{ y }\bws{ x }- ce ( y \bws{ x }+\bws{ y } x )$,
  \item $\sum_{i,j\geq0}a_{ij} \bws{ 1 }^i xy \bws{ 1 }^j$ with the convention that $\bws{ 1 }^0= 1$,
  \item $x \bws{ y }+\bws{ x } y + ax \bws{ 1 } y + bxy,$
  \item $\bws{ x } y + a ( x \bws{ 1 } y - xy \bws{ 1 }),$
  \item $ x \bws{ y }+ a ( x \bws{ 1 } y -\bws{ 1 } xy ) $.
\end{enumerate}
\mlabel{thm:cdto}
\end{theorem}

The differential type operators in terms of convergent rewriting systems and \gsbs were characterizes in~\mcite{GSZ}.
For this, let us recall the monomial order $\ordqd $ in~\mcite{GSZ}.

Let $(Z, \geq)$ be a well-ordered set.  Denote by $\deg_Z(u)$ the
number of $z \in Z$ in $u$ with repetition. Define the order $\ordqd$ on $\opmzm$ as follows. For any $u=u_1\cdots u_m$ and $v=v_1\cdots v_n$, where $u_i$ and $v_j$ are prime.
Define
\begin{equation}
u\ordd v
\mlabel{eq:orddt}
\end{equation}
inductively on $\dep(u)+\dep(v)\geq 0$.
For the initial step of $\dep(u)+\dep(v) = 0$, we have $u,v\in S(Z)$ and define $u\ordd v$ if $u>_{\rm deg-lex} v$, that is
$$(\deg_Z(u), u_1, \ldots, u_m) > (\deg_Z(v), v_1, \ldots, v_n) \, \text{ lexicographically}.$$
For the induction step, if $u = \lc u'\rc$ and $v = \lc v' \rc$, define
$$u\ordd v \,\text{ if }\,  {u'} \ordd {v'}.$$
If $u = \lc u'\rc$ and $v\in Z$, define $u \ordd v$.
Otherwise, define
$$u\ordd v \, \text{ if }\, (\deg_Z(u),  u_1, \ldots, u_m) >(\deg_Z(v),  v_1, \ldots, v_n) \, \text{ lexicographically}.$$

For a system of monic OPIs $\Phi\subseteq\opmx$ and a set $Z$,
denote by $S_\Phi\subseteq \bfk\opm{Z}$ the substitution set:
\begin{equation}
S_\Phi:=\{\phi(u_1,\ldots,u_k)\,|\,\phi\in \Phi,  u_1,\ldots,u_k\in\bfk\opm{Z}\}.
\mlabel{eq:subs}
\end{equation}
If $\Phi = \{\phi\}$ is a singleton set, we write $S_\phi$ for simplicity.

\begin{theorem}\mlabel{thm:drwe}
	\mcite{GSZ}
Let $Z$ be a set and $\phi(x , y ) := \bws{x  y}- N(x  ,y )\in\bfk\opm{\{x,y\}}$ be an OPI satisfying the conditions~\meqref{item:a} and~\meqref{item:b} in Definition~\mref{defn:difftyp}.
With the monomial order $\ordqd$, the following statements are equivalent.
\begin{enumerate}
  \item $\phi(x, y)$ is of differential type OPI.
  \item The term-rewriting system in Eq.~\meqref{eq:difftr0}$:$
\begin{equation*}
\Pia{S_\phi}:=\{q\suba{\bws{u  v}}\rightarrow q\suba{N(u  ,v )}\,|\,q\in\sopma{Z},u,v\in\opm{Z} \}
\end{equation*}
  is convergent. \mlabel{item:difftr}
  \item  The set
$$S_\phi :=\{\phi(u , v ) = \bws{u  v}- N(u  ,v )\,|\, u,v\in \opm{Z}\}$$
is a \gsb in $\bfk\opm{Z}$ with respect to $\ordqd$. \mlabel{item:diffgsb}
\item The free $\phi$-algebra on a set $Z$ is the non-commutative polynomial $\bfk$-algebra $\bfk\langle\Delta(Z)\rangle$ where
$$\Delta(Z):=\{\der{z}{n}\,|\,z\in Z, \,\der{z}{0}:=z,\, \der{z}{n+1}:=\bws{\der{z}{n}},\, n\geq0\},$$
together with the operator $d$ defined by the following recursion:
Let $u = u_1 u_2 \cdots u_k \in M(\Delta(Z))$ for $u_i\in M(\Delta(Z)), 1\leq i\leq k$.
\begin{enumerate}
\item  If $k = 1$, i.e., $u = \der{z}{n}$ for some $i \geq0, z \in Z$, then define $d ( u ) = \der{z}{n+1}$.
\item  If $k \geq 1$, then recursively define $d ( u ) = N ( u_1 , u_2\cdots u_k )$.
\end{enumerate}
\end{enumerate}
\end{theorem}

\subsection{Rota-Baxter type OPIs} \mlabel{ss:rbtype}
Parallel to the case of differential type OPIs, the classification of Rota-Baxter type OPIs was also studied~\mcite{ZGGS}.

\begin{defn}\mcite{ZGGS}
A {\bf Rota-Baxter type OPI} is
$$\phi(x , y ) := \bws{x}\bws{ y}- \bws{B(x , y )}\in\bfk\opm{\{x,y\}},$$
where
\begin{enumerate}
  \item $\bws{B(x, y )}$ is multi-linear in $x$ and $y$ ;
  \mlabel{item:ar}
  \item $\bws{B(x , y )}$ is a normal $\phi$-form, that is, $B(x , y )$  does not contain any subword of the form $\bws{u}\bws{v}$, for any $u,v\in \bfk\opm{\{x,y\}}$;
  \mlabel{item:br}
   \item The term-rewriting system $$\Pia{S_\phi}:=\{q\suba{\bws{u}\bws{ v}}\rightarrow q\suba{\bws{B(u , v )}}\,|\,q\in\sopma{Z},u,v\in\opmzm \}$$ is terminating;
  \item For any $u, v, w \in \opmz$,
${B(B(u,v) , w )} - {B(u , B(v,w) )}\astarrow_{\Pia{S_\phi}} 0.$
  \mlabel{item:cr}
\end{enumerate}
A linear operator satisfying a Rota-Baxter OPI is called a {\bf Rota-Baxter type operator}.
\mlabel{defn:rbtyp}
\end{defn}
The following is a list of Rota-Baxter type operators with operator degrees not exceeding two.
\begin{theorem}\mcite{ZGGS} {\em (Classification of Rota-Baxter type operators)}
For any $c,\lambda\in\bfk$, the OPI $\phi(x,y) := \bws{x}\bws{y} - \bws{B(x , y)}$, where ${B(x , y)}$  is taken from the list below, is of Rota-Baxter type.
\begin{enumerate}
  \item  $x\bws{ y }$ (average operator),
  \item  $\bws{ x } y$ (inverse average operator),
  \item  $ x \bws{ y }+ y \bws{ x },$
  \item  $\bws{ x } y +\bws{ y } x,$
  \item  $x \bws{ y }+\bws{ x } y -\bws{ xy }$ (Nijenhuis operator),
  \item  $x \bws{ y }+\bws{ x } y +\lambda xy$ (Rota-Baxter operator of weight $\lambda$ ),
  \item  $x \bws{ y }- x \bws{ 1 } y +\lambda xy,$
  \item  $\bws{ x } y - x \bws{ 1 } y +\lambda xy,$
  \item  $ x \bws{ y }+\bws{ x } y - x \bws{ 1 } y +\lambda xy$ (generalized Leroux TD operator with weight $\lambda$ ),
  \item  $ x \bws{ y }+\bws{ x } y - xy \bws{ 1 }- x \bws{ 1 } y +\lambda xy, $
  \item  $ x \bws{ y }+\bws{ x } y - x \bws{ 1 } y -\bws{ xy }+\lambda xy,$
  \item  $ x \bws{ y }+\bws{x } y - x \bws{ 1 } y -\bws{ 1 } xy +\lambda xy,$
  \item  $dx \bws{ 1 } y +\lambda xy$ (generalized endomorphisms),
  \item  $dy \bws{ 1 } x +\lambda yx$ (generalized antimorphisms).
\end{enumerate}
 \mlabel{thm:rbtyp}
\end{theorem}
This list of Rota-Baxter type operators may not be complete. See~\mcite{GSZ, ZGGS} for more details.
Similarly, Rota-Baxter type OPIs can be characterized in terms of rewriting systems and \gsbs.

\begin{theorem}\mcite{ZGGS}
Let $Z$ be a set and $\phi(x , y ) := \bws{x}\bws{ y}- \bws{B(x , y)}\in\bfk\opm{\{x,y\}}$ be an OPI satisfying the conditions~\meqref{item:ar} and~\meqref{item:br} in Definition~\mref{defn:rbtyp}.
Let $\geq$ be a monomial order on $\opmzm$ satisfying $\lbar{\phi(u , v )}=\bws{u}\bws{v}$ for $u,v\in\opmzm$. Then the following statements are equivalent.
\begin{enumerate}
  \item $\phi(x, y )$ is of Rota-Baxter type OPI.
  \item The term-rewriting system
  $$\Pia{S_\phi}:=\{q\suba{\bws{u}\bws{ v}}\rightarrow q\suba{\bws{B(u , v )}}\,|\,q\in\sopma{Z},u,v\in\opmzm \}$$
  is convergent.\mlabel{item:rbtr}
  \item  With the monomial order $\geq$, the set
$$S_\phi :=\{\phi(u , v ) = \bws{u }\bws{v}- \bws{B(u , v )}\,|\, u,v\in\opmzm\}$$
is a \gsb in $\opmz$. \mlabel{item:rbgsb}
\end{enumerate}
\mlabel{thm:rbrwe}
\end{theorem}

There is also a uniform construction of free objects for the algebras with Rota-Baxter type operators.

\subsection{The general case}
In the study of \rcp, two developments held the central stage.
The first one is to establish an algebraic framework in which to consider algebraic identities satisfied by a linear operator in \rcp.
This is reviewed in Subsection~\mref{ssec:oalg}.

The first development naturally leads to the question in our understanding of \rcp: what kind of OPIs Rota was interested in?
After the special cases of differential type operators~\mcite{GSZ} and Rota-Baxter operators~\mcite{ZGGS} summarized in the previous two subsections, this question was addressed in terms of rewriting systems and \gsbs in~\mcite{GG}.

\begin{defn}
Let $\geq$ be a monomial order on $\opmzm$ and $S$ be a monic subset of $\opmz$. We define
\begin{equation}
\Pia{S}:=\{q\suba{\lbar{s}}\rightarrow q\suba{R(s)}\,|\, q \in \sopma{Z}, s=\lbar{s}-R(s)\in S\}\subseteq\opmzm\times \opmz.
\mlabel{eq:trs}
\end{equation}
\end{defn}

\begin{defn}
Let $X$ be a set, and let $\Phi\subseteq\opmx$ be a system of monic OPIs. Let $Z$ be a set.
We call $\Phi$ {\bf convergent} if $\Pia{S_\Phi}$ is convergent, where $S_\Phi$ is given in Eq.~(\mref{eq:subs}).
\end{defn}

On the one hand, \rcp can be interpreted in terms of rewriting systems.

\begin{problem}\mcite{GG}
{\em (\rcp via rewriting systems)}. Determine all convergent systems of OPIs.
\mlabel{prob:rpctrs}
\end{problem}

On the other hand, \rcp can also be interpreted in terms of \gsbs.

\begin{defn}
Let $X$ be a set, and let $\Phi\subseteq\opmx$ be a system of monic OPIs. Let $Z$ be a set.
We call $\Phi$ a {\bf \GS system} if ${S_\Phi}$ is a \gsb in $\opmz$.
\end{defn}

\begin{problem}\mcite{GG}
{\em (\rcp via \gsbs)}. Determine all \GS systems of OPIs.
\mlabel{prob:rpcgsb}
\end{problem}

The relationship between reformulations of \rcp is also studied.
For this, we need the relationship between a \gsb of OPIs and a convergent rewriting system of OPIs.

Denote by $\Id(S)$  the ideal of $\opmz$ generated by $S\subseteq\opmz.$ Define
$$\irr(S):=\opmzm\setminus\{q\suba{\lbar{s}}\,|\,q\in\sopma{Z}, s\in S\}.$$

\begin{theorem}\mcite{GG}
Let $Z$ be a set, and let $\geq$ be a monomial order on $\opmzm$. Let $S$
be a monic subset of $\opmz$ and let $\Pia{S}$ be the term-rewriting system from $S$ in Eq.~\meqref{eq:trs}. Then
the following statements are equivalent.
\begin{enumerate}
  \item $\Pia{S}$ is convergent.
  \item $\Pia{S}$ is confluent.
  \item $\Id(S) \cap \bfk\irr(S) = 0$.
  \item $\Id(S) \oplus \bfk\irr(S) =\opmz$.
  \item $S$ is a \gsb in $\opmz$ with respect to $\geq$.
\end{enumerate}
\end{theorem}

\begin{corollary}\mcite{GG}
With a monomial order $\geq$ on $\opmzm$, the two versions Problem~\mref{prob:rpctrs} and
Problem~\mref{prob:rpcgsb} of \rcp are equivalent.
\mlabel{coro:rcpeq}
\end{corollary}

Based on the above corollary, we can conclude that the items~\meqref{item:difftr} and~\meqref{item:diffgsb} in Theorem~\mref{thm:drwe} are equivalent, and the items~\meqref{item:rbtr} and~\meqref{item:rbgsb} in Theorem~\mref{thm:rbrwe} are equivalent.

Beyond the  differential type and Rota-Baxter type OPIs, we have the following OPI.

\begin{lemma}\mcite{GG}
The modified Rota-Baxter OPI
$$\phi(x,y)=\bws{x}\bws{y}-\bws{\bws{x}y}-\bws{x\bws{y}}-\lambda xy \in \bfk\opm{\{x,y\}}, \lambda\in \bfk$$
is  \GS.
\end{lemma}

\begin{problem} \mlabel{pr:mult}
In recent years, associative algebras with multiple linear operators have attracted quite much attention. Here the linear operators are either of different kind or of the same kind and satisfy certain compatibility conditions. Examples of such structures include Rota-Baxter family algebras~\mcite{Guop,ZG} and matching Rota-Baxter algeras~\mcite{GGZ}. Thus it would be interesting to study \rcp for algebras with multiple operators. 	
\end{problem}

\section{\rcp for Lie algebras}
\mlabel{sec:lie}

The \rcp is for linear operators on associative algebras, motivated by the important role played by these operators in mathematics and physics. Historically, a similar role was also played by linear operators on Lie algebras. This motivates adapting the study of \rcp for Lie algebras. Such a study is carried out in~\mcite{ZGG21} utilizing operated Lie algebras developed in~\mcite{QC}. In this section, we give a summary of these developments.

\begin{problem} \mlabel{pr:other}
By taking an approach similar to the one in this section, \rcp should be considered for other algebraic structures, such as commutative algebras, Poisson algebras, dendriform algebras and pre-Lie algebras.
\end{problem}

\subsection{\Olies}
\mlabel{ssec: oplie}
Let us first recall~\mcite{QC} the construction of free \olies (also called Lie $\Omega$-algebras).

\subsubsection{Associative and non-associative Lyndon-Shirshov words }
\mlabel{ssec:aslsw}
Let $X$ be a set. Denote by $S(X)$ the free semigroup on $X$ and by $\nas{X}$ the set of all non-associative (binary) words on $X$.
For any set $Y$, let $\blw{Y}:=\{\blw{y}\,|\,y\in Y\}$ be a disjoint copy of $Y$.
Define
$$\plie{X}_0:=S(X)  \text{  and } \plien{X}_0:=\nas{X}.$$
Assume that we have defined $\plie{X}_{n-1}$ and $\plien{X}_{n-1}$, for any given $n\geq1$, and recursively define
$$\plie{X}_{n}:= S(X\cup \blw{\plie{X}_{n-1}})\,\text{ and }\, \plien{X}_{n}:= \nas{(X\cup \blw{\plien{X}_{n-1}})}.$$
Then
$$\plie{X}_{n-1}\subseteq\plie{X}_{n}\,\text{ and }\, \plien{X}_{n-1}\subseteq\plien{X}_{n}\,\text{ for }\, n\geq1.$$
Define
$$\plie{X}:=\bigcup_{n\geq0}\plie{X}_{n}\,\text{ and }\, \plien{X}:=\bigcup_{n\geq0}\plien{X}_{n}.$$
Elements of $\plie{X}$ (resp. $\plien{X}$) are called the {\bf associative  (resp. non-associative) bracketed words on $X$.}
We represent elements of the $\plie{X}$ (resp. $\plien{X}$) by $w$ (resp. $(w)$).

We collect some basic notations.

\begin{defn}
\begin{enumerate}
\item An element $w$ of $\plie{X}$ can be uniquely expressed in the form
$w=w_1\cdots w_k$ for some $k\geq 1$ and some $w_i\in X\cup \blw{\plie{X}}$, for $ 1 \leq i \leq k$. In this case, we call $w_i$ {\bf prime} and $k$ the {\bf breadth} of $w$, denoted by $|w|$.

\item Define the {\bf depth} of $w\in\plie{X}$ to be ${\rm dep(w)} := {\rm min}\{n\,|\,w\in\plie{X}_n\}.$

\item For any $(w)\in\plien{X}$, there exists a unique $w\in\plie{X}$ by forgetting the brackets of $(w)$. Then we can define the {\bf depth} of $(w)\in\plien{X}$ to be ${\rm dep((w))} := {\rm dep(w)}.$ This agrees with $\min \{n\,|\, (w)\in \plien{X}_n\}$.

\item The {\bf degree} of $w\in\plie{X}$, denoted by ${\rm deg}(w)$, is defined to be the total number of all occurrences of all $x\in X$ and $\blw{~}$ in $w$.
\end{enumerate}
\end{defn}

For example, if $w=\blw{xy\blw{z}y}xy\in\plie{X}$ with $x,y,z\in X$, then
$$|w|=3, \, {\rm dep}(w)=2\,\text{ and }\, {\rm deg}(w)=8.$$
If $(w)=x((\blw{y}(xy))z)\in \plie{X}$ with $x,y,z\in X$, then $$w=x\blw{y}xyz\in \plien{X}\,\text{ and }\, {\rm dep}((w))={\rm dep}(w)=1.$$

It is well known that Lyndon(-Shirshov) words form a linear basis of a free Lie algebra~\mcite{BC,Shi}.
Let $(X, \geq)$ be a well-ordered set. Define lex-order $\geq_{\rm lex}$ on the free monoid $M(X)$ over $X$ by
$1>_{\rm lex} u$ for any nonempty word $u$, and for any $u = xu$ and  $v = yv'$ with $x, y\in X$,
$$ u>_{\rm lex} v\,\text{ if }\,  x > y\, \text{ or }\, x = y\,\text{ and }\, u'>_{\rm lex} v'.$$
\begin{defn} Let $(X, \geq)$ be a well-ordered set.
\begin{enumerate}
\item An associative word $w\in S(X)$ is called an {\bf associative Lyndon-Shirshov word} on $X$ with respect to the lex-order $\geq_{\rm lex}$, if $w =uv >_{\rm lex}vu$ for any decomposition of $w = uv$, where $u,v\in S(X)$.
\item A non-associative word $(w)\in \nas{X}$  is said to be a {\bf non-associative Lyndon-Shirshov word} on $X$ with respect to the lex-order $\geq_{\rm lex}$, if
\begin{enumerate}
  \item $w$ is an associative Lyndon-Shirshov word on $X$;
  \item if $(w) = ((u)(v))$, then both $(u)$ and $(v)$ are non-associative Lyndon-Shirshov words on $X$;
  \item if $(w) = ((u)(v))$ and $(u) = ((u_1 )(u_2 ))$, then $u_2 \leq_{\rm lex} v.$
\end{enumerate}
\end{enumerate}
\end{defn}

Denote by $\alsw{X}$ (resp. $\nlsw{X}$) the set of all associative (resp. non-associative) Lyndon-Shirshov words on $X$ with respect to the lex-order $\geq_{\rm lex}$. For any $w\in\alsw{X}$, there exists a unique
procedure, called {Shirshov-standard bracketing way}~\mcite{BC,Shi}, to give a non-associative Lyndon-Shirshov word $[w]$. Furthermore,
$$\nlsw{X} =\{[w]\,|\,w \in\alsw{X}\}.$$
We alert the reader to the different meanings of the similarly looking notations $\lc w\rc, (w),[w]$. 

We give the following example as an illustration and refer the reader to the original literature for further details.

Take $w=xxyyxy \in \alsw{X}$ on $X = \{x,y\}$ with $x>y$, we have $[w]=((x((xy)y))(xy))\in\nlsw{X}$. Recursively, it is given as follows.
\begin{itemize}
\item  Bracketing the minimal letter $y$ to the previous letters, we obtain a new associative Lyndon-Shirshov word $x(xy)y(xy)$ on $\{x,(xy),y\}$ with $x>_{\rm lex} (xy) >_{\rm lex}y.$
\item  Bracketing the minimal letter $y$ to the previous letters, we get a new associative Lyndon-Shirshov word $x((xy)y)(xy)$ on $\{x,(xy), ((xy)y)\}$ with $x>_{\rm lex}(xy)>_{\rm lex}((xy)y).$
\item  Bracketing the minimal letter $((xy)y)$ to the previous letters, we have a new associative Lyndon-Shirshov word $(x((xy)y))(xy)$ on $\{(x((xy)y)),(xy)\}$ with $(x((xy)y))>_{\rm lex}(xy).$
\item  Bracketing the minimal letter $(xy)$ to the previous letters, we have a new associative Lyndon-Shirshov word $((x((xy)y))(xy))$ on the single point set $\{((x((xy)y))(xy))\}$.
\item  Then $[w]=((x((xy)y))(xy))\in\nlsw{X}$.
\end{itemize}

Let $X$ be a well-ordered set. Let Lie$(X)$ be the Lie subalgebra of commutator Lie algebra $\big(\bfk M(X) , [-,-]:=\mu-\mu\circ \tau\big)$ generated by $X$, where $\bfk  M(X)$ is the free associative algebra on $X$.
It is well-known that Lie$(X)$ is a free Lie algebra on the set $X$ with a linear basis $\nlsw{X}$~\mcite{Reu}.

\subsubsection{Associative and non-associative Lyndon-Shirshov bracketed words }
\mlabel{ssec:anlsw}
In this subsection, we recall the construction of the free \olies in ~\mcite{QC}. In fact, the operated version of non-associative Lyndon-Shirshov words is a linear basis of a free \olie.

\begin{defn}~\mcite{QC}
\begin{enumerate}
\item An {\bf \olie} is a Lie algebra $L$ together with a linear map
$P_L: L\to L$.
\item A {\bf morphism} from an \olie $(L_1,P_{L_1})$ to an \olie $(L_2,P_{L_2})$
is a Lie algebra homomorphism $f : L_1\to L_2$ such that
$f\circ P_{L_1}=P_{L_2}\circ f.$
\end{enumerate}
\end{defn}

The following Deg-lex order $\ordqc$ on $\plie{X}$ is from~\mcite{QC}, which is a monomial order. Let $(X, \geq)$ be a well-ordered set.
Take $u=u_1\cdots u_m$ and $v=v_1\cdots v_n$, where $u_i$ and $v_j$ are prime. Define
\begin{equation}
u\ordc v \mlabel{eq:ordDl}
\end{equation}
inductively on $\dep(u)+\dep(v)\geq 0$.
For the initial step of $\dep(u)+\dep(v) = 0$, we have $u,v\in S(X)$ and define $u\ordc v$ by
$u >_{\rm deg-lex} v$, that is, $$u \ordc v \, \text{ if }\, ({\rm deg}(u), |u|,  u_1, \ldots, u_m) >({\rm deg}(v), |v|,  v_1, \ldots, v_n) \, \text{ lexicographically}.$$
For the induction step, if $u = \lc u'\rc$ and $v = \lc v' \rc$, define
$$u\ordc v \,\text{ if }\, (\deg(u), {u'}) > (\deg(v), {v'}) \, \text{ lexicographically}.$$
Otherwise, define
$$u\ordc v \, \text{ if }\, ({\rm deg}(u), |u|,  u_1, \ldots, u_m) >({\rm deg}(v), |v|,  v_1, \ldots, v_n) \, \text{ lexicographically}.$$

The operated version of Lyndon-Shirshov words was established in~\mcite{QC} with respect to the above specific monomial order $\ordqc$.
Indeed, it can be done with an arbitrary monomial order $\ordq$ on  $\plie{X}$ satisfying,
for prime elements $u_1, \ldots, u_n \in \plie{X}$ and $\sigma\in S_n$,
\begin{equation}
u_1 \cdots u_n  \ordq u_{\sigma(1)} \cdots u_{\sigma(n)} \Longleftrightarrow u_1 \cdots u_n  \succeq_{\rm lex} u_{\sigma(1)} \cdots u_{\sigma(n)},
\mlabel{eq:dlor}
\end{equation}
where $\succeq$ is the restriction of $\ordq$ on the set of prime elements $X\cup\blw{\plie{X}}$.
Eq.~\meqref{eq:dlor} is the key proposition of the monomial order $\ordqc$ used in~\mcite{QC}.
Notice that the order $\ordqc$ is an example of the order $\ordq$.

To formulate \rcp for Lie algebras, we construct free operated Lie algebras. They are defined by a sequence of bracketed words from a recursion. First define
$$\alsbw{X}_0:=\alsw{{X}}$$
and then
$$\nlsbw{X}_0:=\nlsw{{X}}=\{[w]\,|\,w\in \alsbw{X}_0\}$$
with respect to the order $\succeq_{{\rm lex}}$. To complete the recursion, for any given $n\geq 1$, assume that we have defined
$$\alsbw{X}_{n-1}$$
and then
$$\nlsbw{X}_{n-1}=\{[w]\,|\,w\in \alsbw{X}_{n-1}\}$$
with respect to the order $\succeq_{{\rm lex}}$.
Then we first define
$$\alsbw{X}_n:=\alsw{X\cup \blw{\alsbw{X}_{n-1}}}$$
with respect to the order $\succeq_{{\rm lex}}$.
We then denote
$$\big[X\cup \blw{\alsbw{X}_{n-1}}\big]:=\left\{[w]:=
\left\{\left .
\begin{array}{ll}
w, & \hbox{if  $w\in X$}, \\
\blw{[w']}, & \hbox{if $w=\blw{w'}$}
\end{array}
\right.\,\right|\, w \in X\cup \blw{\alsbw{X}_{n-1}}\right\}$$
and define
\begin{equation}
\nlsbw{X}_n:=\nlsw{\big[X\cup \blw{\alsbw{X}_{n-1}}\big]}=\{[w]\,|\,w\in\alsbw{X}_n\}
\mlabel{eq:bnlsbw}
\end{equation}
with respect to the order $\succeq_{\rm lex}$.
We finally denote
$$\alsbw{X}:=\bigcup_{n\geq0}\alsbw{X}_n\,\text{ and }\,\nlsbw{X}:=\bigcup_{n\geq0}\nlsbw{X}_n.$$
Then we have
\begin{equation}
\nlsbw{X}=\{[w]\,|\,w\in\alsbw{X}\}.
\mlabel{eq:astolie}
\end{equation}
Elements of $\alsbw{X}$ (resp. $\nlsbw{X}$) are called the {\bf associative (resp. non-associative) Lyndon-Shirshov bracketed words or Lyndon-Shirshov $\Omega$-words} with respect to $\ordq$.

It follows from Eq.~(\mref{eq:astolie}) that any associative Lyndon-Shirshov bracketed word $w$ is associated to a unique non-associative Lyndon-Shirshov bracketed word $[w]$. For example, for $w=\blw{xyz}\bws{x}\bws{y}\in \alsbwo{X}{\ordq}$ with $x>y>z\in X$, it follows from Eq.~\meqref{eq:bnlsbw} that
      $$[w]=[[\blw{xyz}][\bws{x}][\bws{y}]]=[\blw{[xyz]}\bws{[x]}\bws{[y]}] = [\blw{(x(yz))}\bws{x}\bws{y}]$$
with $\blw{xyz}\succ_{\rm lex}\bws{x}\succ_{\rm lex} \bws{y}.$ For the word $\blw{(x(yz))}\bws{x}\bws{y}$ on
$\{\blw{(x(yz))}, \bws{x}, \bws{y}\}$, appplying the Shirshov-standard bracketing way, we obtain
$[w] = (\blw{(x(yz))}(\bws{x}\bws{y})).$

For $f\in \bfk\alsbw{X}$, define $\lbar{[f]}:=\lbar{f}$ with respect to the order $\ordq$. Notice that $\lbar{[f]}$ needs not be a monomial of $[f]$. For example, if $f=xyz\in \alsbw{X}$ with $x>y>z\in X$, then $[f]=x(yz)$ and $\lbar{[f]}=\lbar{f}=xyz$ is not a monomial of $[f]=x(yz)$.

Denote by $\oplie$ the operated Lie subalgebra of $\bfk\plie{X}$ generated by $X$ under the Lie bracket $[u,v] = uv - vu$.

\begin{lemma} {\em (\mcite{QC, ZGG21})}
The $\oplie$, together with the natural embedding $X\to \oplie$,
is the free \olie on $X$ with a linear basis $\nlsbw{X}$.
\end{lemma}

Based on the above lemma, we propose the following concepts.

\begin{defn}
Let $\phi:=\phi(x_1, \ldots ,x_k)\in\oplie$ with $k\geq1$ and $x_1,\ldots,x_k\in X$. We call $\phi=0$ (or simply $\phi$) an {\bf operated Lie polynomial identity (OLPI)}.
\end{defn}

Using the universal property of the free \olie $\oplie$, for any \olie $(R, P)$ and any map
$$\theta : \{x_1,\ldots ,x_k\}\to R,\quad x_i\mapsto r_i, i=1,\ldots,k,$$
there is a unique morphism
$$\tilde{\theta}:\oplie\to R$$
of \olies that extends the map $\theta$.
By a slight abuse of notation, define the {\bf evaluation map}
$$\phi(r_1, \ldots ,r_k):=\tilde{\theta}(\phi(x_1, \ldots ,x_k)).$$

\begin{defn}
Let $\phi\in \oplie$ and $(R,P)$ be an \olie. We call
$R$ is a {\bf $\phi$-Lie algebra} and $P$ is a {\bf $\phi$-operator}, if
$$\phi(r_1, \ldots ,r_k)=0,\,\forall r_1, \ldots ,r_k\in R.$$
More generally, for a subset $\Phi\subseteq\oplie$, we say that $R$ (resp. $P$) a {\bf $\Phi$-Lie algebra} (resp. {\bf $\Phi$-operator}) if $R$ (resp. $P$) is a $\phi$-Lie algebra (resp. $\phi$-operator) for each $\phi\in\Phi$.
\end{defn}

For example, for $\phi(x,y) = \bws{[[x][y]]} - [\bws{[x]}[y]] - [[x]\bws{[y]}]\in\oplie$,
a $\phi$-algebra is nothing but a differential Lie algebra. 
It can be simplified as $\phi(x,y) =\bws{xy} - \bws{x}y -[x\bws{y}]=\bws{xy} - \bws{x}y -\bws{y}x\in\oplie$. 

\begin{defn}
Let $R$ be a Lie algebra. A linear operator $P$ on $R$ is called {\bf (Lie) algebraic} if there is $0\neq \Phi\subseteq \oplie$ such that $P$ is a $\Phi$-operator.
\mlabel{de:lalgop}
\end{defn}

\subsubsection{Special normal words}

Analogous to $\mathfrak{M}^\star(Z)$, we can define $\sopm{Z}.$
We always use the monomial order $\ordq$ in this subsection.

\begin{lemma}\mcite{QC}
Let $Z$ be a well-ordered set. For any $(u)\in\plien{Z}$, it has a representation
$(u)=\sum_i \alpha_i[u_i]$, where each $\alpha_i\in\bfk$ and $[u_i]\in \nlsbw{Z}$.
In this case, we denote $[(u)]:=\sum_i \alpha_i[u_i]$.
\mlabel{lem:nassb}
\end{lemma}

For example, let $(u)=(\bws{x}\bws{y})\bws{z}$ with $x>y>z$. Then $(u)=\bws{x}(\bws{y}\bws{z})+(\bws{x}\bws{z})\bws{y}$, where $\bws{x}(\bws{y}\bws{z})$ and $(\bws{x}\bws{z})\bws{y}$ are two non-associative Lyndon-Shirshov bracketed words.
\begin{lemma}\mcite{QC,ZGG21}
Let $Z$ be a well-ordered set and let $q\in\sopm{Z}$ and $u$, $q\suba{u}\in\alsbw{Z}$. Then
$[q\suba{u}]=[q'\suba{uc}]$
for some $q'\in\sopm{Z}$ and $c\in\plie{Z}$, where $c$  may be the empty word. Let
\begin{equation}
[q\suba{u}]_u:=[q'\suba{[uc]}]_{[uc]\mapsto [\cdots[[[u][c_1]][c_2]]\cdots[c_m]]}=[q'\suba{[\cdots[[[u][c_1]][c_2]]\cdots[c_m]]}],
\mlabel{eq:snw}
\end{equation}
where $c=c_1c_2\cdots c_m$ with each $c_i\in \alsbw{Z}$ and $c_j \prec_{\rm lex} c_{j+1}$, $1\leq j\leq m-1$. Then
$$[q\suba{u}]_u=q\suba{[u]} + \sum_i \alpha_i q_i\suba{[u]},$$
where each $\alpha_i\in\bfk$, $q_i\in\sopm{Z}$ and $\lbar{[q\suba{u}]_u}=q\suba{u}\ord q_i\suba{u}$.
\end{lemma}

\begin{exam}
Let $Z=\{x,y,z\}$ be a set with $x>y>z$. Let $q=\star z$ and $u=xy$. Then $q\suba{u}=xyz\in\alsbw{Z}$
and so $$[q\suba{u}]=x(yz)=[xyz]=[\star\suba{uz}]$$
is in $\nlsbw{Z}$.
By Eq.~\meqref{eq:snw},
$[q\suba{u}]_u=[[xy][z]]=x(yz)+(xz)y$.
Further by Lemma~\mref{lem:nassb}, $(xy)z=x(yz)+(xz)y$,
and so $[q\suba{u}]_u=x(yz)+(xz)y=(xy)z=q\suba{[u]}$.

\end{exam}

The above result can be extended from bracketed monomial to bracketed polynomials.

\begin{defn}
Let $q\in\sopm{Z}$ and $f\in \opliez\subseteq \bfk\plie{Z}$ be monic.
We call $$[q\suba{{f}}]_{\lbar{f}}:=[q\suba{\lbar{f}}]_{\lbar{f}}\suba{[\lbar{f}]\mapsto f}$$
a {\bf special normal $f$-word} if $q\suba{\lbar{f}}\in \alsbw{Z}$, where
$[q\suba{\lbar{f}}]_{\lbar{f}}$ is defined by Eq.~\meqref{eq:snw}.
\end{defn}

\begin{lemma}
Let $f\in \opliez$ be monic and $q\suba{\lbar{f}}\in \alsbw{Z}$. Then
$$[q\suba{f}]_{\lbar{f}}=q\suba{f}+\sum_i \alpha_i q_i\suba{f},$$
where each $\alpha_i\in\bfk$, $q_i\in\sopm{Z}$ and $q\suba{\lbar{f}}\ord q_i\suba{\lbar{f}}$.
\mlabel{lem:jeq}
\end{lemma}

\subsection{\rcp for Lie algebras}
\mlabel{ssec: rpclie}
With the preparation of the previous subsection, we can formulate \rcp for Lie algebras in the contexts of rewriting systems and \gsbs.

\subsubsection{Rewriting systems}
A family of OLPIs give rise to a rewriting system. Let $\Phi \subseteq \oplie$ be monic.
For each $s\in S_\Phi:=\{\phi(u_1,\ldots,u_k)\,|\,\phi\in \Phi,  u_1,\ldots,u_k\in\bfk\opm{Z}\}$ and $q\in \sopm{Z}$,
if $q\suba{\lbar{s}}\in \alsbw{Z}$, then we may write
$$[q\suba{s}]_{\lbar{s}}=[q\suba{\lbar{s}}]\dps \big(-R([q\suba{s}]_{\lbar{s}})\big),$$
which under the order $\ordq$ can be viewed as a rewriting rule
$$[q\suba{\lbar{s}}]\to R([q\suba{s}]_{\lbar{s}}).$$

Let us give an example for better understanding.

\begin{exam} \mlabel{ex:difflie}
Let $\phi(x,y)=\bws{x}\bws{y}-\bws{\bws{x}y}-\bws{[x\bws{y}]}\in\opliex$. Let $Z$ be a well-ordered set and
$$s=\phi([u], [v]) =\bws{[u]}\bws{[v]}-\bws{\bws{[u]}[v]}-[\bws{{[u]}\bws{[v]}}]$$
with $u\ord v \in\alsbw{Z}$. For $q_1=\star\bws{w}$, $q_2=\star$ and $ {v}\ord {w}\in\alsbw{Z}$, we have $\lbar{s}=\bws{u}\bws{v}$ and
\begin{eqnarray*}
[q_1\suba{s}]_{\lbar{s}}
&=&[q_1\suba{\lbar{s}}]_{\lbar{s}}\suba{[\lbar{s}]\mapsto s}\\
&=&\big[[\bws{[u]}\bws{[v]}][\bws{[w]}]\big]-\big[[\bws{\bws{[u]}{[v]}}][\bws{[w]}]\big]-\big[[\bws{{[u]}\bws{[v]}}][\bws{[w]}]\big]\\
&=&\bws{[u]}\Big(\bws{[v]}\bws{[w]}\Big)+\Big(\bws{[u]}\bws{z}\Big)\bws{[v]}-\bws{[\bws{[u]}{[v]}]}\bws{[w]}-\bws{[{[u]}\bws{[v]}]}\bws{[w]},
\end{eqnarray*}
with $[q_1\suba{\lbar{s}}]=[\bws{[u]}\bws{[v]}\bws{[w]}]=\bws{[u]}(\bws{[v]}\bws{[w]})$. It induces a rewriting rule
$$\bws{[u]}\Big(\bws{[v]}\bws{[w]}\Big)\rightarrow -\Big(\bws{[u]}\bws{z}\Big)\bws{[v]}+\bws{[\bws{[u]}{[v]}]}\bws{[w]} + \bws{[{[u]}\bws{[v]}]}\bws{[w]}.$$
Similarly,
$$[q_2\suba{s}]_{\lbar{s}}=[q_2\suba{\lbar{s}}]_{\lbar{s}}\suba{[\lbar{s}]\mapsto s}=s=\bws{[u]}\bws{[v]}-\bws{\bws{[u]}[v]}-[\bws{{[u]}\bws{[v]}}],$$  with $[q_2\suba{\lbar{s}}]=\bws{[u]}\bws{[v]}$. It induces a rewriting rule
$$\bws{[u]}\bws{[v]}\rightarrow \bws{\bws{[u]}[v]}+[\bws{{[u]}\bws{[v]}}].$$
\end{exam}

\begin{defn}
Let $Z$ be a well-ordered set. Let $S$ be a monic subset of $\opliez$.   With the monomial order $\ordq$ on $\plie{Z}$, we define a term-rewriting system
\begin{equation}
\Pil{S}:=\bigg\{[q\suba{\lbar{s}}]\rightarrow R([q\suba{s}]_{\lbar{s}})\,\bigg|\,s\in S, q\in \sopm{Z},q\suba{\lbar{s}}\in \alsbw{Z}\bigg\}
\subseteq \nlsbw{Z}\times \opliez.
\mlabel{eq:ltrs}
\end{equation}
\end{defn}

\begin{defn}
Let $\Phi\subseteq\oplie$ be a system of monic OLPIs. Let $Z$ be a well-ordered set.
We call $\Phi$ {\bf convergent} if $\Pil{S_\Phi}$ is convergent, where
$$S_\Phi=\{\phi(u_1,\ldots,u_k)\,|\,u_1,\ldots,u_k\in \opliez, \phi\in \Phi\}\subseteq\opliez.$$
\end{defn}

Now we formulate the \rcp for Lie algebras in terms of rewriting systems.

\begin{problem}
{\em (\rcp for Lie algebras via rewriting systems)} Determine all convergent systems of OLPIs.
\mlabel{prob:rpcltrs}
\end{problem}

\subsubsection{\gsbs}
We use the monomial order $\ordq$ in this subsection. In particular, a \gsb in $\opliez$ is with respect to this order.
There are two kinds of compositions.

\begin{defn}
Let $Z$ be a well-ordered set. Let $f, g \in\opliez$ be monic.
\begin{enumerate}
\item  If there are $u, v\in \plie{Z}$ such that $w = \lbar{f}u = v\lbar{g}$ with
$\max\{ |\lbar{f}|, |\lbar{g}|\}< |w| < |\lbar{f}| + |\lbar{g}|$, then
$$\langle f,g \rangle_w:=\langle f,g \rangle^{u,v}_w:= [fu]_{\lbar{f}} - [vg]_{\lbar{g}}$$
is called the {\bf intersection composition of $f$ and $g$ with respect to $w$}.
\mlabel{item:intcompl}
\item  If there is $q\in\sopm{Z}$ such that $w = \lbar{f} = q\suba{\lbar{ g}},$ then
$$\langle f,g \rangle_w:=\langle f,g \rangle^q_w := f - [q\suba{g}]_{\lbar{g}}$$
is called the {\bf including composition of $f$ and $g$ with respect to $w$}.
\mlabel{item:inccompl}
\end{enumerate}
\mlabel{defn:compl}
\end{defn}

\begin{defn}
Let $Z$ be a well-ordered set. Let $S\subseteq\opliez$ be monic.
\begin{enumerate}
\item An element $f\in\opliez$ is called {\bf  trivial modulo $(S, w)$} if
$$f =\sum_i \alpha_i [q_i\suba{s_i}]_{\lbar{s_i}}$$
with $w\ord q_i\suba{\lbar{s_i}}$ for each special normal $s_i$-word $[q_i\suba{s_i}]_{\lbar{s_i}}$, where each $ \alpha_i\in\bfk$, $q_i\in\sopm{Z}$, $s_i\in S$.
In this case, we write $f\equiv 0 \mod(S, w).$
\item We call $S$ a {\bf \gsb} in $\opliez$ with respect to $\ordq$ if, for all pairs $f, g \in S$, every intersection composition of the form $(f,g)^{u,v}_w$ is trivial modulo $(S, w)$, and every including composition of the form $(f, g)^q_w$ is trivial modulo $(S, w)$.
\end{enumerate}
\end{defn}

Denote by $\Idl(S)$  the Lie ideal of $\opliez$ generated by $S\subseteq\opliez.$ Then
\begin{equation*}
\begin{split}
\Idl(S)=~\bigg\{\sum_{i=1}^n \alpha_i[{q_i} \suba{s_i}]_{\lbar{s_i}}\,\bigg|\, n\geq1, \alpha_i\in\bfk, s_i\in S\, \text{ and }\, q_i\in\sopm{Z} \bigg\}.
\end{split}
\end{equation*}
Define
$$\irrl(S):=\big\{[w]\,\big|\, w\in\alsbw{Z}, w\neq q\suba{\lbar{s}}\,\text{ for }\, s\in S\,\text{ and }\, q\in\sopm{Z}\big\}.$$

We have the following Composition-Diamond lemma for \olies.
\begin{lemma}\mcite{QC,ZGG21}{\em (Composition-Diamond lemma for \olies)} Let $Z$ be a well-ordered set. Let $S\subseteq\opliez$ be a monic set.
With the monomial order $\ordq$ on $\plie{Z}$, then the following conditions are equivalent.
\begin{enumerate}
\item $S$ is a  \gsb in $\opliez$ with respect to $\ordq$.
\item  For all $f \neq 0$ in $\Idl ( S )$ , $\lbar{f} = q\suba{\lbar{s}}\in\alsbw{Z}$ for some $q\in\sopm{Z}$ and $s \in $S.
\item  $\opliez=\bfk \irrl( S )\oplus\Idl(S)$ and $\irrl( S )$ is a \bfk-basis of $\opliez/ \Idl(S)$.
\end{enumerate}
\mlabel{lem:cdmol}
\end{lemma}

Then \rcp for Lie algebras can also be interpreted in terms of \gsbs.

\begin{defn}
Let $\Phi\subseteq\oplie$ be a system of monic OLPIs. Let $Z$ be a well-ordered set.
We call $\Phi$ {\bf \GS} if ${S_\Phi}$ is a \gsb in $\opliez$, where
$$S_\Phi=\{\phi(u_1,\ldots,u_k)\,|\,u_1,\ldots,u_k\in \opliez, \phi\in \Phi\}\subseteq\opliez.$$
\end{defn}

\begin{problem}
{\em (\rcp for Lie algebras via \gsbs)} Determine all \GS systems of OLPIs.
\mlabel{prob:rpclgsb}
\end{problem}

\subsubsection{Equivalence of the two formulations}
We will establish a connection between Problem~\mref{prob:rpcltrs} and Problem~\mref{prob:rpclgsb}. See~\mcite{ZGG21} for details.

With the monomial order $\ordq$, we have the following result.

\begin{theorem}
Let $Z$ be a well-ordered set. Let $S$ be a monic set of  $\opliez$, and let $\Pil{S}$ be the term-rewriting system from $S$ in Eq.~\meqref{eq:ltrs}. Then
the following statements are equivalent.
\begin{enumerate}
  \item $\Pil{S}$ is convergent. \mlabel{item,gsbrtsa}
  \item $\Pil{S}$ is confluent. \mlabel{item,gsbrtsb}
  \item $\Idl(S) \oplus \bfk\irrl(S) =\opliez$. \mlabel{item,gsbrtsd}
  \item $S$ is a \gsb in $\opliez$ with respect to $\ordq$.  \mlabel{item,gsbrtse}
\end{enumerate}
\mlabel{lem:liegreq}
\end{theorem}

Now we are ready to give the relationship between the reformulations of \rcp for Lie algebras.
\begin{corollary}
The two versions Problem~\mref{prob:rpcltrs} and
Problem~\mref{prob:rpclgsb} of \rcp for Lie algebras are equivalent.
\mlabel{coro:rcpleq}
\end{corollary}

We now establish relationship between \rcp for associative algebras and Lie algebras.
Chen and Qiu~\mcite{QC} studied the relationship between \gsbs in free associative algebras and \gsbs in
free Lie algebras with respect to the order $\ordqc$.
Indeed this relationship can be generalized to the monomial order $\ordq$.

\begin{lemma}\mcite{QC,ZGG21}
Let $S\subseteq\opliez\subseteq\bfk \plie{Z}$ be monic. With the monomial order $\ordq$ on $\plie{Z}$,
$S$ is \gsb in $\opliez$ if and only if $S$ is \gsb in $\bfk\plie{Z}$.
 \mlabel{lem:gseq}
\end{lemma}

For any element $[w]$ in $\opliez\subseteq\bfk \plie{Z}$,  there is a unique $w$ in $\bfk\plie{Z}$.
Based on this, we have the following equivalence.

\begin{theorem}\mcite{ZGG21}
Let $S\subseteq\opliez\subseteq\opmz$ be monic. With the monomial order $\ordq$ on $\plie{Z}$, the following statements are equivalent.
\begin{enumerate}
\item $\Pil{S}$ is convergent. \mlabel{item:aleqa}
\item $S$ is a \gsb in $\opliez$.\mlabel{item:aleqb}
\item $\Pia{S}$ is convergent.\mlabel{item:aleqc}
\item $S$ is a \gsb in $\opmz$.\mlabel{item:aleqd}
\end{enumerate}
In other words, under the above hypothesis, the Problems~\mref{prob:rpcltrs}, ~\mref{prob:rpclgsb}, ~\mref{prob:rpctrs} and~\mref{prob:rpcgsb} are equivalent to each other with respect to the order $\ordq$.
\mlabel{thm:alrcpeq}
\end{theorem}

\begin{remark}
Let $S\subseteq\opliez\subseteq\bfk \plie{Z}$ be monic OLPIs.
If $S$ give a ``good'' operator on an associative algebra via the classification of Problem~\mref{prob:rpctrs} or Problem~\mref{prob:rpcgsb}, then
$S$ also give a ``good'' operator on a Lie algebra via the classification of Problem~\mref{prob:rpcltrs} or Problem~\mref{prob:rpclgsb}.
\mlabel{coro:rcpasl}
\end{remark}

\subsection{Differential type OLPIs}
We apply Theorem~\mref{thm:alrcpeq} to give some OLPIs that are \gsbs in free \olie.

Parallel to Definition~\mref{defn:difftyp}, we propose the following notion.

\begin{defn}
\mlabel{defn:ldifftyp}
A {\bf differential type OLPI} is
$$\phi(x , y ) := \bws{x y}- [N(x , y )]\in\opliex, $$
 with $x> y$ such that
\begin{enumerate}
\item $[N(x , y )]$ is multi-linear in $x$ and $y$ ;
\mlabel{item:adi}
\item $[N(x , y )]$ is a normal $\phi$-form, that is, $[N(x , y )]$  does not
contain subwords $\bws{[[u][v]]}$ for any $u,v\in \alsbw{\{x,y\}}$ with $u\ord v$;
\mlabel{item:bdi}
\item \mlabel{item:cdi}
For any well-ordered set $Z$ with $u\ord v\ord w \in \alsbw{Z}$,
$$[N([[u][w]],[v])]-
[N([[u][v]],[w])]+ [N([u],[[v][w]])]\astarrow_{\Pil{S_\phi}} 0,$$
where
\begin{equation*}
 \Pil{S_\phi}:=\left\{[q\suba{\bws{u  v}}]\rightarrow R([q\suba{\phi(u , v )}]_{{\bws{u  v}}})\,|\,q\in\sopm{Z}, u, v\in \alsbw{Z}, u\ord v \right\}.
\end{equation*}
\end{enumerate}
A linear operator on a Lie algebra satisfying a differential type OLPI is called a {\bf differential type operator} on a Lie algebra.
\end{defn}

\begin{remark}
The condition~\meqref{item:c} in Definition~\mref{defn:difftyp} is to ensure that $\lc (uv) w\rc = \lc u(vw)\rc$;
while the above condition~\meqref{item:cdi} is for $$\bws{[[u] [v]] [w] } = \bws{ [u] [[v][w] ]} + \bws{[[u][w]] [v]}.$$
\end{remark}

Now we may propose

\begin{problem}{\em (\rcp for Lie algebras: the differential case)} Find all OLPIs of differential
type by finding all expressions $[N ( x , y )]\in \oplie $ of differential type.
\mlabel{pro:dp}
\end{problem}

We now characterize differential type operators in terms of rewriting systems and \gsbs in the free \olies and the free operated associative algebras.

Let $Z$ be a well-ordered set. Define $\der{z}{n}\in\plie{Z}$, $n \geq 0$, recursively by
$$\der{z}{0}:=0,\, \der{z}{n+1}:= \bws{\der{z}{n}}\,\text{ for }\, m\geq0.$$
Denote
\begin{equation}
\Delta(Z):=\{\der{z}{n}\,|\,z\in Z, n\geq0\}.\mlabel{eq:derf}
\end{equation}

Notice that the monomial order $\ordqd$ given in Eq.~\meqref{eq:orddt} satisfies  Eq.~\meqref{eq:dlor}.
\begin{theorem}\mcite{GSZ,ZGG21}
Let $Z$ be a well-ordered set. Let
$$\phi(x , y ) = \bws{x y}- [N(x,y )] \in\opliex \subseteq \bfk\plie{\{x,y\}},$$
 with $x> y$ and $[N(x,y )]$ satisfies the conditions ~\meqref{item:adi} and~\meqref{item:bdi} in Definition~\mref{defn:ldifftyp}.
 With the monomial order $\ordqd$  given in Eq.~\meqref{eq:orddt},
then the following statements are equivalent.
\begin{enumerate}
\item $\phi(x, y)= \bws{x y}- [N(x,y )](\in \bfk\opliex)$ is a differential type OLPI. \mlabel{item:drga}
\item The term-rewriting system induced in Eq.~(\mref{eq:ltrs}):
  $$\Pil{S_\phi}=\left\{[q\suba{\bws{u  v}}]\rightarrow R([q\suba{\phi(u , v )}]_{{\bws{u  v}}})\,|\,q\in\sopm{Z}, u, v\in \alsbwo{Z}{\ordqd}, u\ordd v \right\}$$
  is convergent.\mlabel{item:drgb}

\item $S_\phi=\{\bws{[[u][v]]}-[N([u],[v])]\,|\,u, v\in \alsbwo{Z}{\ordqd}, u\ordd v\}$ is a \gsb in $\opliez$ with respect to the monomial order $\ordqd$. \mlabel{item:drgc}
\item The free $\phi$-Lie algebra on $Z$ is the free Lie algebra $\bfk\nlsbwd{Z}$ on $\Delta(Z)$ together with the operator $d$ on $\bfk\nlsbwd{Z}$.
Here $\Delta(Z)$ is given in Eq.~\meqref{eq:derf}, and $d$ on $\bfk\nlsbwd{Z}$ defined by the following recursion$:$

For any $[u]=[[u_1]\cdots [u_m]]\in \nlsbwd{Z}$ with $u_i\in \Delta(Z)$, $1\leq i \leq n$
\begin{enumerate}
  \item if $m=1$, then $[u]=\der{z}{n}\in \Delta(Z)$ for some $n\geq0$ and $z\in Z$, define $d([u]):=\der{z}{n+1}$;
  \item if $m>1$, then recursively define $d([u]):=[N([u_1],[u_2]\cdots [u_m])]$.
\end{enumerate}\mlabel{item:drgcc}
\item $S_\phi=\{\bws{uv}-N(u,v)\,|\,u,v\in\bfk\plie{Z}\}$ is a \gsb in $\bfk\plie{Z}$ with respect to the order $\ordqd$. \mlabel{item:drgd}
\item The term-rewriting system
  $$\Pia{S_\phi}=\left\{q\suba{\bws{u  v}}\rightarrow q\suba{N(u,v )}\,|\,q\in\sopm{Z}, u, v\in \bfk\plie{Z} \right \}$$
  is convergent.\mlabel{item:drge}
\item $\phi(x, y)= \bws{x y}- N(x,y )(\in \bfk\plie{\{x,y\}})$ is of differential type OPI. \mlabel{item:drgf}
\end{enumerate}
\mlabel{thm:gsbdt}
\end{theorem}

\begin{defn}
The {\bf operator degree} of a monomial in $\oplie$ is the total number that the operator
$\bws{\, }$ appears in the monomial. The {\bf operator degree } of a polynomial $f$ in $\oplie$ is the maximum of the
operator degrees of the monomials appearing in $f$.
\mlabel{defn:olied}
\end{defn}

Next we try to find all differential type OLPIs under a restriction on the number of operators.
Since $1\notin \oplie$, we need to remove the cases for associative algebras involving $1$ in Theorem~\mref{thm:cdto}
and propose the following answer to the Problem~\mref{pro:dp} in the case of operator degree  not exceeding two.

\begin{conjecture}\mcite{ZGG21}
{\em (Classification of differential type OLPIs)}
Let $a , b , c , e\in\bfk$. Every expressions $[N ( x , y )]\in \opliex $ of differential type takes one (or more) of the forms below for $x> y$
\begin{enumerate}
  \item $ b \big( -\bws{ y }x+\bws{ x } y \big)+ c \bws{ x }\bws{ y }+ exy $ where $b^2 = b + ce$,
  \item  $-ce^2 xy + exy -c \bws{ x }\bws{ y }- ce \big( -\bws{ x }y +\bws{ y } x \big)$.
\end{enumerate}
\mlabel{thm:cdtol}
\end{conjecture}
Notice that
\begin{align*}
b \big( -\bws{ y }x +\bws{ x } y \big)+ c \bws{ x }\bws{ y }+ exy  =&\ b \big( [x \bws{ y }]+\bws{ x } y \big)+ c \bws{ x }\bws{ y }+ exy,\\
-ce^2 xy + exy -c \bws{ x }\bws{ y }- ce \big( -\bws{ x }y +\bws{ y } x \big)
=&\ ce^2 [yx] + exy + c [\bws{ y }\bws{ x }]- ce \big( [y \bws{ x }]+\bws{ y } x \big).
\end{align*}\

By Theorems~\mref{thm:cdto} and~\mref{thm:gsbdt} and the fact that $1\notin\oplie$, we have

\begin{corollary}
Let $[{N(x , y)}]\in \opliex \subseteq \bfk\plie{\{x,y\}}$ be from the list in Conjecture~\mref{thm:cdtol}. Then all the statements in Theorem~\mref{thm:gsbdt} hold.
\end{corollary}

\subsection{Rota-Baxter type OLPIs}
Next we consider Rota-Baxter type operator identities for Lie algebras.

\begin{defn}
\mlabel{defn:lrbtyp}
A {\bf Rota-Baxter type OLPI} is
a $$\phi(x , y ) := \bws{x}\bws{ y}- \bws{[B(x , y )]}\in\opliex$$
with $x> y$ such that
\begin{enumerate}
  \item ${[B(x , y )]}$ is multi-linear in $x$ and $y$ ;
  \mlabel{item:arrb}
  \item ${[B(x , y )]}$ is a normal $\phi$-form, that is, $[B(x , y )]$  does not contain the subword  $[\bws{[u]}\bws{[v]}]$ for any $u,v\in \alsbw{X}$ with $u\ord v$;
  \mlabel{item:brrb}
   \item The term-rewriting system
   $$\Pil{S_\phi}=\big\{[q\suba{\bws{u}\bws{  v}}]\rightarrow [R([q\suba{\phi(u , v )}]_{{\bws{u}\bws{  v}}})]\,|\,q\in\sopm{Z}, u, v\in \alsbwo{Z}{\ordq}, u\ord v \big \}$$ is terminating; \mlabel{item:drrb}
  \item   \mlabel{item:crrb}
 For any well-ordered set $Z$ with $u\ord v\ord w \in \alsbw{Z}$,
$$[B([B([u],[w])],[v])]-[B([B([u],[v])],[w])]+ [B([u],[B([v],[w])])]\astarrow_{\Pil{S_\phi}} 0,$$
where 
$$\Pil{S_\phi}:=\big\{[q\suba{\bws{u}\bws{  v}}]\rightarrow R([q\suba{\phi(u , v )}]_{{\bws{u}\bws{  v}}})\,|\,q\in\sopm{Z}, u, v\in \alsbwo{Z}{\ordq}, u\ord v \big \}.$$
\end{enumerate}
A linear operator on a Lie algebra that satisfies a Rota-Baxter type OLPI is called a {\bf Rota-Baxter type operator} on a Lie algebra.
\end{defn}

\begin{remark}
The condition~\meqref{item:cr} in Definition~\mref{defn:rbtyp} is to insure $(\bws{u}\bws{v}) \bws{w} = \bws{u}(\bws{v} \bws{w})$;
while the above condition~\meqref{item:crrb} is for $$[\bws{[u]}\bws{ [v]}] \bws{[w] } = \bws{ [u]} [\bws{[v]}\bws{[w] }] + [\bws{[u]}\bws{[w]}] \bws{[v]}.$$
\end{remark}

\begin{problem} {\em (\rcp for Lie algebras: the Rota-Baxter case)} Find all OLPIs of Rota-Baxter
type by finding all expressions ${[B ( x , y )]}\in \opliex $ of Rota-Baxter type.
\mlabel{pro:rbp}
\end{problem}

We give some criteria for Rota-Baxter type OLPIs.

\begin{theorem}\mcite{ZGG21, ZGGS} Let $Z$ be a well-ordered set.
Let
$$\phi(x , y ) = \bws{x}\bws{ y}- \bws{[B(x , y )]}\in \opliex \subseteq \bfk\plie{\{x,y\}},$$
with $x> y$ and $\phi(x,y)$ satisfies the conditions~\meqref{item:arrb} and~\meqref{item:brrb} in Definition~\mref{defn:lrbtyp}.
Let $\ordq$ be the monomial order defined in Eq.~\meqref{eq:dlor} satisfying $\lbar{\phi(u , v)}= \bws{u}\bws{ v}$ for $u\ord v\in \alsbwo{Z}{\ordq}$.
Then the following statements are equivalent.
\begin{enumerate}
  \item $\phi(x, y) = \bws{x}\bws{ y}- \bws{[B(x , y )]}(\in \opliex)$ is of Rota-Baxter type OLPI.\mlabel{item:rdrga}
  \item The term-rewriting system
  $$\Pil{S_\phi}=\big\{[q\suba{\bws{u}\bws{  v}}]\rightarrow R([q\suba{\phi(u , v )}]_{{\bws{u}\bws{  v}}})\,|\,q\in\sopm{Z}, u, v\in \alsbwo{Z}{\ordq}, u\ord v \big \}$$
  is convergent.\mlabel{item:rdrgb}
  \item $S_\phi=\{ [\bws{[u]}\bws{ [v]}]- \bws{[B([u] , [v] )]}\,|\,u, v\in \alsbwo{Z}{\ordqd}, u\ordd v\}$ is a \gsb in $\opliez$ with respect to the order $\ordq$.\mlabel{item:rdrgc}
  \item The set
  $$\irrl(S):=\big\{[w]\,\big|\, w\in\alsbwo{Z}{\ordq}, w\neq q\suba{{{\bws{u}\bws{  v}}}}\,\text{ for }\, q\in\sopm{Z}, u, v\in \alsbwo{Z}{\ordq}, u\ord v\big\}$$
  is a linear basis of the free $\phi$-Lie algebra $\opliez/ \Idl(S).$\mlabel{item:rdrgcc}
  \item $S_\phi=\{  \bws{u}\bws{ v}- \bws{B(u , v )}\,|\,u,v\in\bfk\plie{Z} \}$ is a \gsb in $\bfk\plie{Z}$ with respect to the order $\ordq$. \mlabel{item:rdrgd}
   \item The term-rewriting system
  $$\Pia{S_\phi}=\big\{q\suba{\bws{u  v}}\rightarrow q\suba{\bws{B(u,v )}}\,|\,q\in\sopm{Z}, u, v\in \bfk\plie{Z} \big \}$$
  is convergent.\mlabel{item:rdrge}
   \item $\phi(x, y) = \bws{x}\bws{ y}- \bws{B(x , y )}(\in \bfk\plie{\{x,y\}})$ is of Rota-Baxter type OPI. 
\mlabel{item:rdrgf}
\end{enumerate}
\mlabel{thm:gsbrbt}
\end{theorem}

Again since $1\notin \oplie$, we omit the types involving $1$ in Theorem~\mref{thm:rbtyp}
and propose the following answer to Problem~\mref{pro:rbp} in the case of operator degree not exceeding two.

\begin{conjecture} \mcite{ZGG21} {\em (Classification of Rota-Baxter type operators on Lie algebras)}
For any $\lambda\in\bfk$, the expressions $[{B(x , y)}]\in \opliex$ in the list below are of Rota-Baxter type, for $x> y$
\begin{enumerate}
  \item  $-\bws{y}x$ (average operator),
  \item  $\bws{ x } y$ (inverse average operator),
  \item  $\bws{ x } y+\bws{ y } x$,
  \item  $-\bws{y}x-\bws{ x }y $,
  \item  $-\bws{y}x+\bws{ x } y -\bws{ xy }$ $($Nijenhuis operator$)$,
  \item  $\bws{ x } y -\bws{y}x +\lambda xy$ $($Rota-Baxter operator of weight $\lambda$$)$.
\end{enumerate}
\mlabel{conj:rbt}
\end{conjecture}
Notice that
\begin{align*}
-\bws{y}x=&\ [x\bws{y}],\\
-\bws{y}x-\bws{ x }y=&\ [x\bws{y}]+[y\bws{x}],\\
-\bws{y}x+\bws{ x } y -\bws{ xy }=&\ [x\bws{y}]+\bws{x}y-\bws{xy},\\
\bws{ x } y-\bws{y}x+\lambda xy=&\ \bws{ x } y+[x\bws{y}]+\lambda xy.
\end{align*}

The next corollary follows from Theorems~\mref{thm:rbtyp} and~\mref{thm:gsbrbt} and the fact that $1$ is not in $\oplie$.

\begin{corollary}
Let $[{B(x , y)}]\in \opliex \subseteq \bfk\plie{\{x,y\}}$ be from the list in Conjecture~\mref{conj:rbt}. Then all the statements in Theorem~\mref{thm:gsbrbt} hold.
\end{corollary}

\subsection{Modified Rota-Baxter OLPI}
Recall the monomial order $\ordqc$ on $\plie{Z}$ defined in Eq.~\meqref{eq:ordDl}.
The modified Rota-Baxter OLPI of weight $\lambda$ is
$$\phi(x,y)= \bws{x}\bws{ y}-[\bws{x \bws{ y }}]-[\bws{\bws{ x } y}]-\lambda xy \in\opliex,$$
with $x> y$.
When $\lambda=-1$, the modified Rota-Baxter Lie algebra  is a Lie algebra $L$ equipped with a linear map $P : L \to L$ satisfying
$$[P(x)P(y)]=P([P(x)y])+P([xP(y)])-[xy]\,\text{ for }\,x,y\in L.$$
The term modified Rota-Baxter operator was adapted from the modified Yang-Baxter equation in Lie algebra~\mcite{STS}. See~\mcite{GLS,ZGG1,ZGG2} for more recent developments on this operator.
Note the subtle difference between this operator and the Rota-Baxter operator.

\begin{theorem}\mcite{GG,QC,ZGG21}
 Let $$\phi(x,y)= \bws{x}\bws{ y}-[\bws{x \bws{ y }}]-[\bws{\bws{ x } y}]-\lambda xy \in\opliex$$
be an OLPI with $x> y$.
Let $Z$ be a well-ordered set. With the order $\ordqc$ on $\plie{Z}$, the following statements are equivalent.
\begin{enumerate}
  \item $S_\phi=\{\bws{u}\bws{ v}-\bws{u \bws{ v }}-\bws{\bws{ u } v}-\lambda uv\,|\,u,v\in\plie{Z}\}$ is a \gsb in $\bfk\plie{Z}$ with respect to the monomial order $\ordqc$.
  \item
  The term-rewriting system $$\Pia{S_\phi}=\big\{q\suba{\bws{u}\bws{  v}}\rightarrow q\suba{\bws{\bws{u}  v}+\bws{{u}\bws{  v}}+\lambda uv}\,|\,q\in\sopm{Z}, u, v\in \plie{Z}\big\}$$
  is convergent.
  \item $S_\phi=\{[\bws{[u]}\bws{ [u]}]-[\bws{[u] \bws{ [v] }}]-[\bws{\bws{ [u] } [v]}]-\lambda [[u][v]] \,|\,\,u, v\in \alsbwo{Z}{\ordqc}, u\ordc v\}$ is a \gsb in $\opliez$ with respect to the monomial order $\ordqc$.
  \item the term-rewriting system $$\Pil{S_\phi}=\big\{[q\suba{\bws{u}\bws{  v}}]\rightarrow R([q\suba{\phi(u,v)}]_{{\bws{u}\bws{v}}})\,|\,q\in\sopm{Z}, u, v\in \alsbwo{Z}{\ordqc}, u\ordc v \big \}$$
  is convergent.
\end{enumerate}
\end{theorem}

\section{Appendix: Early work on \rcp}
\mlabel{sec:free}
For completeness, we summarize the early work of Freeman~\mcite{Fr} on \rcp. The goal of his work is different from ours in that he focused on a classification of several operators known to the author at that time, rather than attempting a more complete classification of operator identities. By considering isomorphisms of direct product algebras corresponding to operator identities,
	the classification scheme of several operator identities was obtained as follows~\mcite{Fr}.
	\begin{description}
		\item[$\mathscr{H}$] Identities of homomorphism type:
		\begin{align}
			\text{(homomorphism operator)\quad\quad} P(x)P(y)&=P(xy), \tag{$ \cdot $}\mlabel{eq:tphom}\\
			\text{(RBO of weight -1)\quad\quad} P(x)P(y)&=P(P(x)y+xP(y)- xy),\tag{$\circ$}\mlabel{eq:tprb1}\\
			\text{(Hilbert transform)\quad\quad} P(x)P(y)&=P(P(x)y+xP(y))+xy.  \tag{$\odot$}\mlabel{eq:tpmrb}
		\end{align}
		\mlabel{item:tph}
		\item[$\mathscr{D}$] Identities of derivation type:
		\begin{align}
			\text{(derivation)\quad\quad} P(xy)&=P(x)y+xP(y), \tag{$ \bigtriangleup $}\mlabel{eq:tpder}\\
			\text{(RBO of weight 0)\quad\quad} P(x)P(y)&=P(P(x)y+xP(y)),\tag{$\bigtriangledown$}\mlabel{eq:tprb0}\\
			\text{(Reynolds operators)\quad\quad} P(x)P(y)&=P(P(x)y+xP(y)-P(x)P(y)).  \tag{$\square$}\mlabel{eq:tprn}
		\end{align}
		\mlabel{eq:tpd}
		\item[$\mathscr{A}$] Right Averaging  Identity:
		\begin{equation}
			\text{(right averaging)\quad\quad} P(x)P(y)=P(xP(y)).\tag{$\times$}
			\mlabel{eq:tpra}
		\end{equation}
	\end{description}
The first two of the homomorphism type are of the Rota-Baxter type in our sense (See Subsection~\ref{ss:rbtype}); the third one is not, but is obtained from a Rota-Baxter operator after a linear transformation and is called the modified Rota-Baxter operator in recent literature as noted in the introduction.
	
It was noted in~\mcite{Fr} that an operator $P: A\to A$ is a homomorphism operator on an algebra $A$ if and only if the graph of $P$:
	$${\rm G(P)} = \{(a, P(a))\,|\, a\in A\}$$
	is a subalgebra of algebra $A\times A$,
	where the multiplication on $A\times A$ is given by $$(a,b)\cdot (a',b'):=(aa',bb').$$
	There are analogous interpretations for the other operators defined by Eqs.~\eqref{eq:tprb1}, ~\eqref{eq:tpmrb}, ~\eqref{eq:tpder},  ~\eqref{eq:tprb0}, ~\eqref{eq:tprn} and ~\eqref{eq:tpra}, where the multiplications on $A\times A$ are given respectively by
	\begin{equation*}
		\begin{split}
			(a,b)\circ (a',b'):=&~(ab'  +  ba' - aa',bb'),\\
			(a,b)\odot(a',b'):=&~(ab'+  ba', bb' - aa'),\\
			(a,b)\bigtriangleup (a',b'):=&~(aa', ab'  +  ba'),\\
			(a,b)\bigtriangledown (a',b'):=&~(ab'  +  ba', bb'),\\
			(a,b)\square (a',b'):=&~(ab'  +  ba' - bb', bb'),\\
			(a,b)\times (a',b'):=&~(ab', bb').\\
		\end{split}
	\end{equation*}
	Notice that all these multiplications are associative on $A\times A$.
	
	Based on the close relationshiop between the operator $P:A \to A$ and the subalgebra
	$G(P)$ of $A\times A$, Freeman posed the following concept to classify these operators.
	\begin{defn}~\mcite{Fr}
		Let $P: A\to A$ be an operator on an algebra $A$, and let $(A\times A, \ast)$ be an algebra.
		The operator $P$ is called an {\bf $\ast$-operator} if the graph $G(P)$ is a subalgebra of $(A\times A, \ast)$.
	\end{defn}
	In terms of $\ast$-operators and isomorphisms, Freeman gave a standard to classify operators.
	\vskip0.1in
	
	\noindent{\bf Freeman's standard:} Let $P:A\to A$ (resp. $P':A'\to A'$) be an $\ast$-operator (resp. $\ast'$-operator) on
	the algebra $A$ (resp. $A'$). If there is an isomorphism $(A\times A, \ast) \cong (A\times A, \ast')$,
	then the $\ast$-operator $P$ and the $\ast'$-operator $P'$ are said to be in the same class.
	\vskip0.1in
	
	The isomorphism in Freeman's standard can be achieved by the following concept of a $\sigma$-transform.
	
	\begin{defn}~\mcite{Fr}
		Let $A$ and $(A\times A, \ast)$ be algebras, and let $\sigma$ be a linear automorphism on $A\times A$.
		Define a multiplication $\ast_\sigma$ on $A\times A$ by
		\begin{equation}
			(a,b)\ast_\sigma (a',b'):=\sigma\Big(\sigma^{-1}(a,b)\ast \sigma^{-1}(a',b')\Big).
			\mlabel{eq:smat}
		\end{equation}
		Then $(A\times A, \ast_\sigma)$ is an algebra, called the $\sigma$-{\bf transform} of $(A\times A, \ast)$.
	\end{defn}
	
	Let us emphasize that the defining equation~(\mref{eq:smat}) is equivalent to that $\sigma:(A\times A, \ast)\to (A\times A, \ast_\sigma)$ is an isomorphism.
	
	Thanks to $\sigma$-transforms, Freeman showed that the homomorphism operator, RBO of weight -1 and Hiblert transform are in the same type.
	
	\begin{theorem}~\mcite{Fr}
		Let $u$ and $v $ be two non-zero divisors of $A$.
		Let $\sigma_1$ and $\sigma_2$ be two linear automorphisms on $A\times A$ with
		$$\sigma_1(a,b)=(b, b-u a v )\,\text{ and }\, \sigma_2(a,b)=(b+u a v , b-u a v ).$$
		Then the algebra  $(A\times A, \cdot)$ is
		the $\sigma_1$-transform of $(A\times A, \circ)$ and the $\sigma_2$-transform of $(A\times A, \odot)$,
		where
		\begin{align*}
			(a,b)\circ (a',b')&=(av  b'v ^{-1}  +  u^{-1}bu a' - av  u a',bb'),\\
			(a,b)\odot(a',b')&=(av  b'v ^{-1}+  u^{-1}bu a', bb' +u av u a'v ).
		\end{align*}
	\end{theorem}
	
	Further Freeman obtained that the derivation, RBO of weight 0 and Reynolds operator are in the same type.
	
	\begin{theorem}~\mcite{Fr}
		Let $\sigma_1$ and $\sigma_2$ be two automorphism on $A\times A$ with
		$$\sigma_1(a,b)=(b, a)\,\text{ and }\, \sigma_2(a,b)=(b, b-a).$$
		Then the algebra  $(A\times A, \bigtriangleup)$ is the $\sigma_1$-transform of $(A\times A, \bigtriangledown)$ and the $\sigma_2$-transform of $(A\times A, \square)$.
	\end{theorem}

The approach of Freeman summarized here made critical use of the associative products on a direct product. This is quite similar in style to the formulations of the differential type and Rota-Baxter type operators which also made critical use of the associativity. This observation naturally leads to the following questions.
\begin{problem}
\begin{enumerate}
\item Can the approach of Freeman be applied to the classification of the other known linear operators, in particular the differential type and Rota-Baxter type operators? For instance, does the Nijenhuis operator also defines an associative operation on the direct product $A\times A$?
\item In the other direction, can the idea of Freeman be applied in the general approach of \rcp?
\item Does the approach of Freeman work for linear operators on Lie algebras? From the close analogy of the associative algebra case and the Lie algebra case in the previous two sections, the answer is likely to be affirmative.
\end{enumerate}
\end{problem}

\smallskip
\noindent
{\bf Acknowledgments.}
This work is supported by the National Natural Science Foundation of
China (Grant No. 11771190, 1861051 and 12071191), and the Natural Science Foundation of Gansu
Province (Grant No. 20JR5RA249) and the Natural Science Foundation of Shandong Province
(ZR2020MA002).

\end{document}